\documentclass[11pt,reqno]{amsart}
\usepackage{amsmath}
\usepackage{amssymb}
\usepackage{amsfonts}

\def\squarebox#1{\hbox to #1{\hfill\vbox to #1{\vfill}}}

\newcommand{\Z}{{\mathbb Z}}

\newcommand{\R}{{\mathbb R}}
\newcommand{\C}{{\mathbb C}}
\newcommand{\N}{{\mathbb N}}

\renewcommand{\Re}{\mathop{\rm Re}\nolimits}
\renewcommand{\Im}{\mathop{\rm Im}\nolimits}

\theoremstyle{plain}

\newtheorem{thm}{Theorem}

\newtheorem{lem}{Lemma}
  
\newtheorem{prop}{Proposition}

\newtheorem{deff}{Definition}

\begin{document}

\def\sn{{\bf S}^{n-1}}
\def\ts{\tilde{\sigma}}
\def\ss{{\mathcal S}}
\def\aa{{\mathcal A}}
\def\cc{{\mathcal C}}
\def\tpp{\widetilde{{\mathcal P}}}
\def\iu{\underline{i}}
\def\lc{{\mathcal L}}
\def\pxi{\phi + (\xi + i u)\psi_p}
\def\pp{\mathcal P}
\def\rr{\mathcal R}
\def\hU{\widehat{U}}
\def\mt{\Lambda}
\def\hLa{\widehat{\mt}}
\def\ep{\epsilon}
\def\tPhi{\widetilde{\Phi}}
\def\hR{\widehat{R}}
\def\oV{\overline{V}}
\def\uu{\mathcal U}
\def\tz{\tilde{z}}
\def\hz{\hat{z}}
\def\hd{\hat{\delta}}
\def\ty{\tilde{y}}
\def\hs{\hat{s}}
\def\hc{\hat{\cc}}
\def\hC{\widehat{C}}
\def\hh{\mathcal H}
\def\tf{\tilde{f}}
\def\of{\overline{f}}
\def\trr{\tilde{r}}
\def\tr{\tilde{r}}
\def\tts{\tilde{\sigma}}
\def\tVl{\widetilde{V}^{(\ell)}}
\def\tVj{\widetilde{V}^{(j)}}
\def\tVo{\widetilde{V}^{(1)}}
\def\tVj{\widetilde{V}^{(j)}}
\def\tPsi{\tilde{\Psi}}
 \def\tp{\tilde{p}}
 \def\tVjo{\widetilde{V}^{(j_0)}}
\def\tvj{\tilde{v}^{(j)}}
\def\tVjj{\widetilde{V}^{(j+1)}}
\def\tvl{\tilde{v}^{(\ell)}}
\def\tVt{\widetilde{V}^{(2)}}
\def\Lo{\; \stackrel{\circ}{L}}
\def\tg{\tilde{g}}

\def\ii{{\imath }}
\def\jj{{\jmath }}
\vspace*{0,8cm}
\def\saa{\Sigma_A^+}
\def\sAA{\Sigma_{\aa}^+}
\def\Lip{\mbox{\rm Lip}}
\def\clip{C^{\mbox{\footnotesize \rm Lip}}}
\def\lip{\mbox{{\footnotesize\rm Lip}}}
\def\Vol{\mbox{\rm Vol}}

\def\ccm{\cc^{(m)}}
\def\oo{\mbox{\rm O}}
\def\ooo{\oo^{(1)}}
\def\oot{\oo^{(2)}}
\def\ooj{\oo^{(j)}}
\def\fo{f^{(1)}}
\def\ft{f^{(2)}}
\def\fj{f^{(j)}}
\def\wo{w^{(1)}}
\def\wt{w^{(2)}}
\def\wj{w^{(j)}}
\def\Vo{V^{(1)}}
\def\Vt{V^{(2)}}
\def\Vj{V^{(j)}}

\def\Uo{U^{(1)}}
\def\Ut{U^{(2)}}
\def\Ul{U^{(\ell)}}
\def\Uj{U^{(j)}}
\def\wl{w^{(\ell)}}
\def\Vl{V^{(\ell)}}
\def\Ujj{U^{(j+1)}}
\def\wjj{w^{(j+1)}}
\def\Vjj{V^{(j+1)}}
\def\Ujo{U^{(j_0)}}
\def\wjo{w^{(j_0)}}
\def\Vjo{V^{(j_0)}}
\def\vj{v^{(j)}}
\def\vl{v^{(\ell)}}

\def\gl{\gamma_\ell}
\def\id{\mbox{\rm id}}
\def\piU{\pi^{(U)}}
\def\bs{\bigskip}
\def\ms{\medskip}
\def\Int{\mbox{\rm Int}}
\def\diam{\mbox{\rm diam}}
\def\di{\displaystyle}
\def\dist{\mbox{\rm dist}}
\def\ff{{\cal F}}
\def\i{{\bf i}}
\def\pr{\mbox{\rm pr}}
\def\co{\; \stackrel{\circ}{C}}
\def\la{\langle}
\def\ra{\rangle}
\def\supp{\mbox{\rm supp}}
\def\Arg{\mbox{\rm Arg}}
\def\Int{\mbox{\rm Int}}
\def\II{{\mathcal I}}
\def\e{\emptyset}
\def\endofproof{{\rule{6pt}{6pt}}}
\def\con{\mbox{\rm const }}
\def\Box{\spadesuit}
\def\be{\begin{equation}}
\def\ee{\end{equation}}
\def\beqn{\begin{eqnarray*}}
\def\eeqn{\end{eqnarray*}}
\def\MM{{\mathcal M}}
\def\tmu {\tilde{\mu}}
\def\Pr{\mbox{\rm Pr}}
\def\Prf{\mbox{\footnotesize\rm Pr}}
\def\htau{\hat{\tau}}
\def\btau{\overline{\tau}}
\def\hr{\hat{r}}
\def\tF{\widetilde{F}}
\def\tG{\widetilde{G}}
\def\trho{\tilde{\rho}}

\title[Two parameters]{Spectral estimates for Ruelle transfer operators with two parameters and applications}
\author[V. Petkov]{Vesselin Petkov}
\address{Universit\'e de Bordeaux, Institut de Math\'ematiques de Bordeaux, 351,
Cours de la Lib\'eration, 
33405  Talence, France}
\email{petkov@math.u-bordeaux1.fr}
\author[L. Stoyanov]{Luchezar Stoyanov}
\address{University of Western Australia, School of Mathematics 
and Statistics, Perth, WA 6009, Australia}
\email{luchezar.stoyanov@uwa.edu.au}
\thanks{VP was partially supported by ANR project Nosevol BS01019 01}

\def\ts{\tilde{\sigma}}
\def\ss{{\mathcal S}}
\def\aa{{\mathcal A}}
\def\R{{\mathbb R}}
\def\C{{\mathbb C}}
\def\iu{\underline{i}}
\vspace*{0,8cm}
\def\saa{\Sigma_A^+}
\def\sAA{\Sigma_{\aa}^+}
\def\lc{{\mathcal L}}
\def\pxi{\phi + (\xi + i u)\psi_p}
\def\sa{\Sigma_A}
\def\ssa{\Sigma_A^+}
\def\ssn{\sum_{\sigma^n x = x}}
\def\lc{{\mathcal L}}
\def\lo{{\mathcal O}}
\def\oo{{\mathcal O}}
\def\mt{\Lambda}
\def\ep{\epsilon}
\def\ms{\medskip}
\def\bs{\bigskip}
\def\diam{\mbox{\rm diam}}
\def\rr{\mathcal R}
\def\pp{\mathcal P}
\def\hU{\widehat{U}}
\def\hz{\hat{z}}
\def\tz{\tilde{z}}
\def\ty{\tilde{y}}
\def\i{{\bf i}}
\def\pp{{\mathcal P}}
\def\ff{{\mathcal F}_{\theta}}
\def\tG{\tilde{G}}
\def\hw{\hat{w}}
\def\ep{\epsilon}
\def\rt{R^{\tau}}
\maketitle

\begin{abstract} For $C^2$ weak mixing Axiom A flow $\phi_t: M \longrightarrow M$ on a Riemannian manifold $M$ and a basic set 
$\Lambda$ for $\phi_t$ we consider the Ruelle transfer operator $L_{f - s \tau + z g}$, where $f$ and $g$ are real-valued 
H\"older functions on $\Lambda$, $\tau$ is the roof function and $s, z \in \C$ are complex parameters. Under some assumptions 
about $\phi_t$ we establish estimates for the iterations of this Ruelle operator in the spirit of the estimates for operators 
with one complex parameter (see \cite{D}, \cite{St2}, \cite{St3}). Two cases are covered: (i) for arbitrary H\"older $f,g$
when $|\Im z| \leq B |\Im s|^\mu$ for some constants $B > 0$, $0 < \mu < 1$ ($\mu = 1$ for Lipschitz $f,g$), (ii) for Lipschitz
$f,g$ when $|\Im s| \leq B_1 |\Im z|$ for some constant $B > 0$ . Applying these estimates,  we obtain a non zero analytic
extension of the zeta function $\zeta(s, z)$ for $P_f - \epsilon < \Re (s) < P_f$ and $|z|$ small enough with simple pole 
at $s = s(z)$. Two other applications are considered as well: the first concerns the Hannay-Ozorio de Almeida sum formula, 
while the second deals with the asymptotic of the counting function $\pi_F(T)$ for weighted primitive periods of the flow $\phi_t.$

\end{abstract}

\section{Introduction}
\renewcommand{\theequation}{\arabic{section}.\arabic{equation}}
\setcounter{equation}{0}

Let $M$ be a $C^2$ complete (not necessarily compact) 
Riemannian manifold,   and let $\phi_t: M \longrightarrow M, \: t \in \R,$ 
be a $C^2$ weak mixing Axiom A flow (see \cite{B2}, \cite{PP}). Let $\mt$ be a {\it basic set} for $\phi_t$,  i.e. $\mt$ is a compact
$\phi_t-$ invariant subset of $M$, $\phi_t$ is hyperbolic and transitive on $\mt$ and $\mt$ is locally maximal, i.e. there 
exists an open neighborhood $V$ of $\mt$ in $M$ such that $\mt = \cap_{t\in \R} \phi_t(V)$. The restriction of the flow $\phi_t$ on 
$\Lambda$ is a hyperbolic flow \cite{PP}. For any $x \in M$ let $W_{\epsilon}^s(x), W_{\epsilon}^u(x)$ be the local stable and
unstable manifolds through $x$,  respectively (see \cite{B2}, \cite{KH},  \cite{PP}). 

When $M$ is compact and $M$ itself is a basic set, $\phi_t$ is called an {\it Anosov flow}.
It follows from the hyperbolicity of $\mt$  that if  $\epsilon_0 > 0$ is sufficiently small,
there exists $\ep_1 > 0$ such that if $x,y\in \mt$ and $d (x,y) < \ep_1$, 
then $W^s_{\ep_0}(x)$ and $\phi_{[-\ep_0,\ep_0]}(W^u_{\ep_0}(y))$ intersect at exactly 
one point $[x,y ] \in \mt$  (cf. \cite{KH}). This means that there exists a unique 
$t\in [-\ep_0, \ep_0]$ such that $\phi_t([x,y]) \in W^u_{\ep_0}(y)$. Setting $\Delta(x,y) = t$, 
defines the so called {\it temporal distance function}.

In the paper we will use the set-up and some arguments from \cite{St2}.  First, as in \cite{St2}, we fix a 
{\it (pseudo-) Markov partition} $\rr = \{ R_i\}_{i=1}^k$ of {\it pseudo-rectangles} 
$R_i = [U_i  , S_i ] =  \{ [x,y] : x\in U_i, y\in S_i\}$. Set $R = \cup_{i=1}^{k} R_i,\: U = \cup_{i=1}^k U_i$. Consider the
{\it Poincar\'e map} $\pp: R \longrightarrow R$, defined by  $\pp(x) = \phi_{\tau(x)}(x) \in R$, where
$\tau(x) > 0$ is the smallest positive time with $\phi_{\tau(x)}(x) \in R$. 
The function $\tau$  is the so called {\it first return time}  associated with $\rr$. 
Let  $\sigma : U   \longrightarrow U$ be the {\it shift map} given by
$\sigma  = \piU \circ \pp$, where $\piU : R \longrightarrow U$ is the {\it projection} along stable leaves. 
Let $\hU$ be the set of those points $x\in U$ such that $\pp^m(x)$ is not a boundary point of a rectangle
for any integer $m$. In a similar way define $\hR$.
Clearly in general $\tau$ is not continuous on $U$, however under the assumption that
the holonomy maps are Lipschitz (see Sect. 3) $\tau$ is {\it essentially Lipschitz} on $U$ 
in the sense that there exists a constant $L > 0$ such that if $x,y \in U_i \cap \sigma^{-1}(U_j)$ 
for some $i,j$, then $|\tau(x) - \tau(y)| \leq L\, d(x,y)$.
The same applies to $\sigma : U \longrightarrow U$.


The hyperbolicity of the flow on $\mt$ implies the existence of constants $c_0 \in (0,1]$ and $\gamma_1 > \gamma_0 > 1$ such that
\begin{equation} \label{eq:1.1}
 c_0 \gamma_0^m\; d (u_1,u_2) \leq  d (\sigma^m(u_1), \sigma^m(u_2)) \leq \frac{\gamma_1^m}{c_0} d (u_1,u_2)
\end{equation}
whenever $\sigma^j(u_1)$ and $\sigma^j(u_2)$ belong to the same  $U_{i_j}$  for all $j = 0,1 \ldots,m$.

Define a $k \times k$ matrix $A = \{A(i, j)\}_{i, j = 1}^k$ by
$$A(i, j) = \begin{cases} 1 \:\:{\rm if}\: \pp({\rm Int}\:R_i)\cap {\rm Int}\: R_j \neq \emptyset,\\
0 \:\:{\rm otherwise}. \end{cases}$$
It is possible to construct a Markov partition ${\mathcal R}$ so that $A$ is irreducible and aperiodic (see \cite{B2}).\\
 Introduce 
$R^{\tau} = \{(x, t) \in R \times \R:\: 0 \leq t \leq \tau(x)\}/ \sim,$ where by $\sim$ we identify the points $(x, \tau(x))$ and $(\sigma x, 0).$ 
One defines the suspended flow  $\sigma_t^{\tau}(x, s) = (x, s + t)$  on $R^{\tau}$ taking into account the identification $\sim.$
 For a H\"older continuous function $f$ on $R$, the pressure $\Pr(f)$ with respect to $\sigma$ is defined as
$$\Pr(f) = \sup_{m \in {\mathcal M}_{\sigma}} \big\{ h(\sigma, m) + \int f d m\big \},$$
where ${\mathcal M}_{\sigma}$ denotes the space of all $\sigma$-invariant Borel probability measures and $h(\sigma, m)$ is 
the  entropy of $\sigma$ with respect to $m$.  We say that $f$ and $g$ are cohomologous and we denote this by 
$f \sim g$ if there exists a continuous function $w$ such that $f = g + w \circ \sigma - w.$ For a function $v$ on $R$ one defines 
$$v^n(x) : = v(x) + v(\sigma(x)) +...+ v(\sigma^{n-1}(x)).$$

Let $\gamma$ denote a primitive periodic orbit of $\phi_t$ and let $\lambda(\gamma)$ denote its least period. Given a H\"older 
function $F : \Lambda \longrightarrow \R$, introduce the weighted period 
$\lambda_F(\gamma) = \int_0^{\lambda(\gamma)} F(\phi_t(x_{\gamma})) dt,$ where $x_{\gamma} \in \gamma.$ Consider the weighted
 version of the dynamical zeta function (see Section 9 in \cite{PP})
$$\zeta_{\phi}(s, F): = \prod_{\gamma}\Bigl( 1 - e^{\lambda_F(\gamma) - s\lambda(\gamma)}\Bigr)^{-1}.$$
 Denote by $\pi(x, t): R^{\tau} \longrightarrow \Lambda$ the semi-conjugacy projection
which is one-to-one on a residual set  and $\pi(t, x) \circ \sigma_t ^\tau  = \phi_t \circ \pi(t, x)$ (see \cite{B2}). Then following the results in \cite{B2}, \cite{BR}, a closed  $\sigma$-orbit 
$\{x, \sigma x, ..., \sigma^{n-1}x\}$  is projected to a closed orbit $\gamma$ in $\Lambda$ with a least period
$$\lambda(\gamma) = \tau^n(x): = \tau(x) + \tau(\sigma(x)) +...+ \tau(\sigma^{n-1}(x)).$$

Passing to the symbolic model $R$ (see \cite{B2}, \cite{PP}), the analysis of $\zeta_{\varphi}(s, F)$ is reduced to that of the
Dirichlet series
$$\eta(s) = \sum_{n=1}^{\infty} \frac{1}{n} \sum_{\sigma^n x = x}  e^{f^n(x) - s \tau^n(x)}.$$
with a H\"older continuous function $f(x) = \int_0^{\tau(x)} F(\pi(x, t)) dt: R \longrightarrow \R.$
On the other hand, to deal with certain problems (see Chapter 9 in \cite{PP} and \cite{PoS2}) it is necessary to study a more general series
$$\eta_g(s) =\sum_{n=1}^{\infty} \frac{1}{n} \sum_{\sigma^n x = x} g^n(x) e^{f^n(x) - s \tau^n(x)}$$
with a H\"older continuous function $G: \Lambda \longrightarrow \R$ and 
$g(x) = \int_0^{\tau(x)} G(\pi(x, t)) dt: R \longrightarrow \R$. 
For this purpose it is convenient to examine the zeta function

\begin{equation} \label{eq:1.2}
\zeta(s, z): = \prod_{\gamma} \Bigl(1 - e^{\lambda_F(\gamma) - s \lambda(\gamma) + z \lambda_G(\gamma)}\Bigr)^{-1} 
= \exp\Bigl( \sum_{n=1}^{\infty} \frac{1}{n} \sum_{\sigma^n x = x}  e^{f^n(x) - s \tau^n(x)+ z g^n(x)}\Bigr)
\end{equation}
depending on two complex variables $s, z \in \C.$ Formally, we get
$$\eta_g(s) =  \frac{\partial \log \zeta(s, z)}{\partial z} \Bigl\vert_{z = 0}.$$

{\bf Example 1.} If $G = 0$ we obtain the classical Ruelle {\it dynamical zeta function}
$$\zeta_{\phi}(s) = \prod_{\gamma}\Bigl( 1 - e^{- s\lambda(\gamma)}\Bigr)^{-1}.$$
Then $Pr(0) = h,$ where $h > 0$ is the topological entropy of $\phi_t$ and $\zeta_{\phi}(s)$ is absolutely convergent for $\Re s > h$ (see Chapter 6 in \cite{PP}).\\

{\bf Example 2.} Consider the expansion function $E: \Lambda \longrightarrow \R$ defined by
$$E(x) : = \lim _{t \to 0} \frac{1}{t} \log | {\rm Jac}\:(D\phi_t \vert_{E^u(x)})|,$$
where the tangent space $T_x(M)$ is decomposed as $T_x(M) = E^s(x) \oplus E^0(x) \oplus E^u(x)$  with $E^s(x), E^u(x)$ tangent to stable and instable manifolds through $x$, respectively.
Introduce the function $\lambda^u(\gamma) = \lambda_E(\gamma)$
 and  define $f: R \longrightarrow \R$ by
$$f(x) = - \int_0^{\tau(x)} E(\pi(x, t)) dt.$$
Then we  have $-\lambda^u(\gamma) = f^n(x)$ , $f$ is H\"older continuous function and $Pr(f) = 0$ (see \cite{BR}). Consequently, the series
\begin{equation} \label{eq:1.3}
\sum_{n=1}^{\infty} \frac{1}{n} \sum_{\sigma^n x = x}  e^{f^n(x) - s \tau^n(x)}
\end{equation}
is absolutely convergent for $\Re s > 0$ and nowhere zero and analytic for $\Re s \geq 0$ except for a simple pole at 
$\Re s = 0$ (see Theorem 9.2 in \cite{PP}). The roof functions $\tau(x)$ is constant on stable leaves of rectangles $R_i$ of the Markov family ${\mathcal R}$, 
so we can assume that $\tau(x)$ depends only on $x \in U.$ By a standard argument (see \cite{PP}) we can replace $f$ in (\ref{eq:1.3}) by a H\"older function $\hat{f}(x)$ which depends only 
on $x \in U$ so that $f \sim \hat{f}$.  Thus the series (\ref{eq:1.3}) can be written by functions $\hat{f},\: \tau$ depending only on $x \in U$. We keep the notation $f$ below 
assuming that $f$ depends only on $x \in U.$ The analysis of the analytic continuation of (\ref{eq:1.3}) is based on spectral  estimates for the iterations of the Ruelle operator
$$L_{f - s \tau} v(x) = \sum_{\sigma y = x} e^{f(y) - s \tau(y)} v(y),\: v \in C^{\alpha}(U),\:s \in \C.$$
(see for more details \cite{D}, \cite{PoS1}, \cite{St2}, \cite{St3}, \cite{W}).\\

{\bf Example 3.} Let $f, \tau$ be real-valued H\"older functions and let $P_f > 0$ be the unique real number such that 
$Pr(f - P_f \tau) = 0.$ Let $g(x) = \int_0^{\tau(x)} G(\pi(x, t))dt,$ where $G : \Lambda \longrightarrow \R$ is a H\"older 
function. Then if the suspended flow $\sigma^\tau_t$ is weak-mixing, the function (\ref{eq:1.2}) is nowhere zero analytic 
function for $\Re s > P_f$ and $z$ in a neighborhood of 0 (depending on $s$) with nowhere zero analytic extension to 
$\Re s = P_f \:(s \neq P_f)$ for small $|z|$. This statement is just Theorem 6.4 in \cite{PP}.  To examine the analytic 
continuation of $\zeta(s, z)$ for $P_f - \eta_0 \leq \Re s, \: \eta_0 > 0$ and small $|z|$, it is necessary to establish and to exploit some spectral estimates 
for the iterations of the Ruelle operator
\begin{equation} \label{eq:1.4}
L_{f - s \tau + z g}v(x) = \sum_{\sigma y = x} e^{f(y) - s \tau(y) + z g(y)} v(y),\: v \in C^{\alpha}(U), \: s \in \C, z \in \C.
\end{equation}
The analytic continuation of $\zeta(s, z)$ for small $|z|$ and that of $\eta_g(s)$ play a crucial role in the argument in 
\cite{PoS2} concerning the Hannay-Ozorio de Almeida sum formula for the geodesic flow on compact negatively curved surfaces. 
We deal with the same question for Axiom A flows on basic sets in Sect. 7.\\

{\bf Example 4.} In the paper \cite{KS} the authors  examine for Anosov flows the spectral properties of the Ruelle operator (\ref{eq:1.4}) 
with $f = 0$ and $z = \i w, \: w \in \R,$ as well as  the analyticity of the corresponding L-function $L(s, z).$ The properties of the  
Ruelle operator
$$L^n_{f - (P_f + a + \i b)\tau + \i w},\:w \in \R,\: n \in \N, $$
are also rather important in the paper \cite{Wa} dealing with the large deviations for Anosov flows. Here as above $P_f \in \R$
is such that $Pr(f - P_f \tau) = 0.$ However, it is important to note that in \cite{KS} and \cite{Wa} the analysis of the Ruelle
operators covers mainly the domain $\Re s \geq P_f$ and there are no results treating the spectral properties for 
$P_f - \eta_0 \leq \Re s < P_f$ and $ z = \i w, \: w \in \R.$ To our best knowledge the analytic continuation of the function 
$\zeta(s, z)$ for these values of $s$ and $z$ has not been investigated in the literature so far which makes it quite difficult to 
obtain sharper results.\\

In this paper under some hypothesis on the flow $\phi_t$ (see Sect. 3 for our standing assumptions) we prove spectral 
estimates for the iterations of the Ruelle operator $L_{f - s \tau + z g}^n$ with {\bf two complex parameters} $s, z \in \C$. 
These estimates are in the spirit of those obtained in \cite{D}, \cite{St1}, \cite{St2}, \cite{St3} for the Ruelle operators with 
{\bf one complex parameter} $s \in \C$. On the other hand, in this analysis some new difficulties appear when  
$|\Im s| \to \infty$ and $|\Im z| \to \infty.$ First we prove in Theorem 5 spectral estimates in the case of arbitrary H\"older
continuous functions $f,g,$ when there exist constants $B > 0$ and $0 < \mu < 1$ such that $|\Im z| \leq B |\Im s|^\mu$ and
$|\Im s| \geq b_0 > 0$. When $f,g$ are Lipschitz one can take $\mu = 1$. This covers completely the case when $|z|$ is bounded 
and the estimates have the same form as those for operators with one complex parameter. Moreover, these estimates are sufficient for 
the applications in \cite{PP} and \cite{PoS2} when $|z|$ runs in a small neighborhood of 0 (see Sect. 6 and 7).  In Sect. 5 we
deal with the case when $f,g$ are Lipschitz and there exists a constant $B_1 > 0$ 
such that $|\Im s| \leq B_1 |\Im z|$ (see Theorem 6).\\

To study the  analytic continuation of $\zeta(s, z)$ for $P_f - \eta_0 < \Re s < P_f$, we need a generalization of the so called
Ruelle's lemma which yields a link between the convergence by packets of a Dirichlet series like (\ref{eq:1.3})  and the estimates of the iterations of the corresponding Ruelle operator. The reader may consult \cite{W} 
for the precise result in this direction and the previous works (\cite{R}, \cite{PoS1}, \cite{N}), treating this question. 
For our needs in this paper we prove in Sect. 2 an analogue of Ruelle's lemma for Dirichlet series with two complex parameters 
following the approach in \cite{W}. Combining Theorem 4 with the estimates in Theorem 5 (b), we obtain the following

\begin{thm} Assume  the standing assumptions in Sect. $3$ fulfilled for a basic set $\Lambda.$ Then for any  H\"older continuous functions
$F, G: \Lambda \longrightarrow \R$ there exists $\eta_0 > 0$ such that the function $\zeta(s, z)$ admits a non zero analytic
continuation for 
$$(s, z) \in \{(s, z) \in \C^2:\: P_f - \eta_0 \leq \Re s,\: s \neq s(z),\: |z| \leq \eta_0\}$$
with a simple pole at $s(z).$ The pole $s(z)$ is determined as the root of the equation $Pr(f - s \tau + z g) = 0$ with respect to $s$ for $|z| \leq \eta_0.$
\end{thm}
 


Applying the results of Sects. 4, 5, we  study also the analytic continuation of $\zeta(s, \i w)$ for  $P_f - \eta_0 < \Re s $ and 
$ w \in \R$, $|w| \geq \eta_0$, in the case when $F, G: \Lambda \longrightarrow \R$ are Lipschitz functions (see Theorem 7). 
This analytic continuation combined with the arguments in \cite{Wa} opens some new
perspectives for the investigation of sharp large deviations for Anosov flows with exponentially shrinking intervals in the spirit 
of \cite{PeS}.\\

Our first application concerns the so called Hannay-Ozorio de Almeida sum formula (see \cite{HO}, \cite{Pa}, \cite{PoS3}). 
Let $\phi_t: M \longrightarrow M$ be the geodesic flow on the unit-tangent bundle over a compact negatively curved surface $M$. 
In \cite{PoS3} it was proved that there exists $\epsilon > 0$ such that if $\delta(T) = {\mathcal O} (e^{-\epsilon T})$, 
for every H\"older continuous function $G: M \longrightarrow \R$, we have  
\begin{equation} \label{eq:1.5}
\lim_{T \to +\infty} \frac{1}{\delta(T)} \sum_{T - \frac{\delta(T)}{2} \leq \lambda(\gamma) 
\leq T + \frac{\delta(T)}{2}}   \lambda_G(\gamma) e^{-\lambda^u(\gamma)} = \int G d \mu,
\end{equation}
where  the notations $\lambda(\gamma), \lambda_G(\gamma)$ and $\lambda^u(\gamma)$ for a primitive periodic orbit $\gamma$ are introduced above,
while $\mu$ is the unique $\phi_t$-invariant probability measure which is absolutely continuous with respect to the volume 
measure on $M$. The measure $\mu$  is called SRB (Sinai-Ruelle-Bowen) measure (see \cite{BR}). Notice that in the above case 
the Anosov flow $\phi_t$ is weak mixing and $M$ is an attractor. Applying Theorem 1 and the arguments in \cite{PoS3}, we prove 
the following

\begin{thm}  Let $\Lambda$ be an attractor, that is there exists an open neighborhood $V$ of $\Lambda$ such that 
$\Lambda =\cap_{t \geq 0} \phi_t(V).$ Assume the standing assumptions of Sect. $3$ fulfilled for the basic set $\Lambda$. 
Then there exists $\epsilon > 0$ such that if $\delta(T) = {\mathcal O} (e^{-\epsilon T})$, then for every H\"older function 
$G: \Lambda \longrightarrow \R$ the formula $(\ref{eq:1.5})$ holds with the SRB measure $\mu$ for $\phi_t.$
\end{thm}

Our second application concerns the counting function
$$\pi_F(T) = \sum_{\lambda(\gamma) \leq T} e^{\lambda_F(\gamma)},$$
where $\gamma$ is a primitive period orbit for $\phi_t: \Lambda \longrightarrow \Lambda$, $\lambda(\gamma)$ is the least period 
and $\lambda_F(\gamma) = \int_0^{\lambda(\gamma)} F(\phi_t(x_{\gamma})) dt,\: x_{\gamma} \in \gamma.$ For $F = 0$ we obtain 
the counting function $\pi_0(T) = \#\{\gamma: \: \lambda(\gamma) \leq T\}.$ These counting functions have been studied in many works 
(see \cite{PoS1} for references concerning $\pi_0(T)$ and \cite{PP}, \cite{Po1} for the function $\pi_F(T)$). The study
of $\pi_F(T)$ is based on the analytic continuation of the function
$$\zeta_F(s) = \prod_{\gamma} \Bigl( 1 - e^{\lambda_F(\gamma) - s \lambda(\gamma)}\Bigr)^{-1}, \: s \in \C$$
which is just the function $\zeta(s, 0)$ defined above. We prove the following

\begin{thm} Let $\Lambda$ be a basic set and let $F: \Lambda \longrightarrow \R$ be a H\"older function. Assume the standing
assumptions of Sect. $3$ fulfilled for $\Lambda$. Then there exists $\epsilon > 0$ such that
$$\pi_F(T) = li(e^{Pr(F) T})(1 + {\mathcal O} (e^{-\epsilon T})),\: T \to \infty ,$$
where $li(x):= \int_2^{x} \frac{1}{\log y} dy \sim \frac{x}{\log x},\: x \to +\infty.$
\end{thm}
In the case when $\phi_t: T^1(M) \longrightarrow T^1(M)$ is the geodesic flow on the unit tangent bundle $T^1(M)$ of a compact 
$C^2$ manifold $M$ with negative section curvatures which are $\frac{1}{4}$-pinching the above result has been established in
\cite{Po1}. Following \cite{St2}, \cite{St3}, one deduces that the special case of a geodesic flow in \cite{Po1} is covered by  Theorem 3.

\section{Ruelle lemma with two complex parameters}
\renewcommand{\theequation}{\arabic{section}.\arabic{equation}}
\setcounter{equation}{0}

Let $B(\hU)$ be the {\it space of  bounded functions} $q : \hU \longrightarrow \C$ with its standard norm  
$\|q\|_0 = \sup_{x\in \hU} |g(x)|$. Given a function $q \in B(\hU)$, the  {\it Ruelle transfer operator } 
$L_q : B(\hU) \longrightarrow B(\hU)$ is defined by $\di (L_qh)(u) = \sum_{\sigma(v) = u} e^{q(v)} h(v)\;.$
If $q \in B(\hU)$ is Lipschitz on $\hU$ with respect to the Riemann metric, then  $L_q$ preserves the space $\clip(\hU)$
of {\it Lipschitz functions} $q: \hU \longrightarrow \C$. Similarly, if $q$ is $\nu$-H\"older for some
$\nu > 0$, the operator $L_q$ preserves the {\it space $C^\nu(\hU)$ of $\nu$-H\"older functions} on $\hU$.
In this section we assume that $g, \tau$ and $f$ are real-valued $\nu-$H\"older continuous functions on $\hat{U}$. 
Then we can extend these functions as H\"older continuous on $U$. 

 We define the Ruelle operator $L_{g - sr + z f}: C^\nu(\hat{U}) \longrightarrow C^\nu(\hat{U})$ by
$$L_{f - s\tau + zg} v(x) = \sum_{\sigma y = x} e^{f(y) - s \tau(y) + zg(y)} v(y),\: s, z \in \C.$$
Next, for $\nu > 0$  define the $\nu$-norm on a set $B \subset U$ by
$$|w|_{\nu} = \sup\Bigl\{ \frac{|w(x) - w(y)|}{d(x, y)^{\nu}}: x, y \in B \cap U_i,\: i = 1,...,k,\: x \neq y\Bigr\} .$$
Let 
$$\|w\|_{\nu} = \|w\|_{\infty} + |w|_{\nu} ,$$ 
and denote by $\|. \|_{\nu}$ be the corresponding norm for operators. Let $\chi_i(x)$ be the characteristic function of $U_i.$  

Introduce the sum 
$$Z_n(f - sr + zg): = \sum_{\sigma^n x = x} e^{f^n(x) - s \tau^n (x) + z g^n(x)}.$$

Our purpose is to prove the following statement which can be considered as Ruelle's lemma with two complex parameters.

\begin{thm} For every Markov leaf $U_i$ fix an arbitrary point $ x_i \in U_i$. Then for every $\epsilon > 0$ and sufficiently 
small $a_0 > 0,  c_0 > 0$ there exists a constant $C_{\epsilon} > 0$ such that
\begin{eqnarray} \label{eq:2.1}
\Bigl| Z_n(f - s\tau + z g) - \sum_{i=1}^k L_{f - s  \tau + zg}^n \chi_i(x_i)\Bigr| \nonumber \nonumber \\
\leq C_{\epsilon}(1 + |s|)(1 + |z|) \sum_{m = 2}^n \|L_{f- s\tau + zg}^{n-m}\|_{\nu} \gamma_0^{-m \nu} 
e^{m(\epsilon + Pr (f-a\tau + c g))},\: \forall n \in \N
\end{eqnarray}
for $s = a + \i b,\: z = c + \i w ,\: |a| \leq a_0, |c| \leq c_0.$
\end{thm}

The proof of this theorem follows  that of Theorem 3.1 in \cite{W} with some modifications. We have to take into account the presence 
of a second complex parameter $z$. Given a string $\alpha = (\alpha_0,...,\alpha_{n-1})$ of symbols $\alpha_j$ taking the values 
in  $\{1,...,k\}$, we say that $\alpha$ is an admissible word if $A(\alpha_j, \alpha_{j+1}) = 1$ for all $0 \leq j \leq n-1$. Set 
$|\alpha| = n$ and define the cylinder of length $n$ in the leaf $U_{\alpha_0}$ by
$$U_{\alpha} = U_{\alpha_0} \cap \sigma^{-1} U_{\alpha_1} \cap...\cap \sigma^{-(n-1)}U_{\alpha_{n-1}}.$$
Each $U_i$ is a cylinder of length 1. Next we introduce some other words (see Section in \cite{W}). Given a word 
$\alpha = (\alpha_0,...,\alpha_{n-1})$ and $i = 1,...,k$, if $A(\alpha_{n-1}, i) = 1$ and $A(i, \alpha_0) = 1$, we define
$$\alpha i = (\alpha_0,..., \alpha_{n-1}, i),\:  i \alpha = (i, \alpha_0,..., \alpha_{n-1}), \: 
\bar{\alpha} = (\alpha_0,...,\alpha_{n-2}).$$
We have the following 

\begin{lem} Let $w$ be a $\nu$-H\"older real-valued function on $U$. Let $x$ and $y$ be on the same cylinder $U_{\alpha}$ with 
$|\alpha| = m.$ Then there exists a constant $B > 0$ depending only on $w, \nu$ and the constants $c_0$ and $\gamma_0$ in $(\ref{eq:1.1})$ such that
$$|w^m(x) - w^m(y)| \leq B (d(\sigma^{m-1}x, \sigma^{m-1} y))^{\nu}.$$
\end{lem}

The proof is a repetition of that of Lemma 2.5 in \cite{W} and we leave the details to the reader.

\begin{prop} Let $m \geq 1$ and let $w$ be a function which is $\nu$-H\"older continuous on all cylinder of length $m + 1$. 
Then for the transfer operator $L_{f - s\tau + zg}$ we have
$$L_{f - s\tau + zg}: = \oplus_{|\alpha| = m+1}C^{\nu}(U_{\alpha}) \ni w  \longrightarrow L_{f - s\tau + z g} w \in 
\oplus_{|\alpha| = m} C^{\nu}(U_{\alpha}).$$
\end{prop}

{\bf Proof.} Let $w$ be $\nu$-H\"older on $U_{i\alpha}$ for all $i$ such that $A(i, \alpha_0) = 1$. Let 
$x, y \in {\rm Int}\: U_{\alpha}$ and let $|U| = \max_{i=1,,,k} \diam (U_i)$. Then
$$|L_{f - s\tau + zg}w(x) - L_{f - s\tau + zg}w(y)| $$
$$= \Bigl|\sum_{A(i, \alpha_0)= 1} e^{f(ix)-s\tau(ix)+zg(ix)} w(ix) 
- \sum_{A(i, \alpha_0)= 1} e^{f(iy)-s\tau(iy)+zg(iy)} w(iy)\Bigr|$$
$$\leq \sum_{A(i, \alpha_0) = 1} |e^{-s\tau(iy)}|\Bigr( | e^{s\tau(iy) - s\tau(ix)} -1| |e^{f(iy) + zg(iy)} w(ix)|
 + |e^{f(iy) + zg(iy)}w(iy) - e^{f(ix) + zg(ix)} w(ix)|\Bigr)$$
$$\leq e^{a_0 |\tau|_{\infty}}\sum_{A(i,\alpha_0) = 1} \Bigl(|s| |\tau|_{\beta} e^{a_0|\tau|_{\nu} 
|U|^{\nu}}e^{|f|_{\infty} +c_0|g|_{\infty}}|w|_{\infty} + |e^{f(iy) + zg(iy)}w(iy) - e^{f(ix) + zg(ix)} w(ix)|\Bigr).$$
Repeating this argument, we get
$$\sum_{A(i,\alpha_0) = 1} |e^{f(iy) + zg(iy)}w(iy) - e^{f(ix) + zg(ix)} w(ix)|$$
$$\leq e^{c_0|g|_{\infty}}\sum_{A(i, \alpha_0) = 1} \Bigl( |z| |g|_{\nu} e^{c_0|g|_{\nu} |U|^{\nu}} e^{|f|_{\infty}} |w|_{\infty} + |e^{f(iy)} w(iy) - e^{f(ix)} w(ix)|\Bigr)$$
and we conclude that
$$|L_{f - s\tau + zg}w(x) - L_{f - s\tau + zg}w(y)| \leq C |w|_{\nu} d(x, y)^{\nu}.$$
\hfill\qed

\medskip

\def\xa{x_{\alpha}}

Now, as in \cite{W}, we will choose in every cylinder $U_{\alpha}$ a point $\xa \in U_{\alpha}.$ For the reader's 
convenience we recall the choice of $\xa.$ \\
$\:\:$
(1) If $U_{\alpha}$ has an $n$-periodic point, then we take $\xa \in U_{\alpha}$ so that $\sigma^n \xa = \xa.$\\
(2) If $U_{\alpha}$ has no $n$-periodic point and $n > 1$ we choose $\xa \in U_{\alpha}$ arbitrary so that  $\xa \notin \sigma(U_{\alpha_{n-1}}).$\\
(3) if $|\alpha| = n = 1$, then we take $\xa = x_i$, where $i = \alpha_0$ and $x_i \in U_i$ is one of the points fixed in Theorem 4.\\

Let $\chi_{\alpha}$ be the characteristic function of $U_{\alpha}.$ Then  Lemma 3.4 and Lemma 3.5 in \cite{W} are applied without any change and we get
$$Z_n(f - s \tau + z g) = \sum_{|\alpha| = n} (L_{f - s\tau + zg}^n \chi_{\alpha})(\xa).$$

\begin{prop} We have
\begin{eqnarray} \label{eq:2.2}
Z_n(f - s\tau + zg) - \sum_{i=1}^k L_{f - s\tau + zg}^n \chi_i(x_i) \nonumber \\
= \sum_{m = 2}^n\Bigl( \sum_{|\alpha| = m} L_{f - s\tau + zg}^n \chi_{\alpha}(x_{\alpha}) - \sum_{|\beta| 
= m -1} L^n_{f - s\tau + zg} \chi_{\beta} (x_{\beta})\Bigl).\end{eqnarray}
\end{prop}
The proof is elementary by using the fact that 
$$\sum_{i=1}^k (L^n_{f - s\tau + zg} \chi_{U_i}) ( x_i) = \sum_{|\alpha| = 1} (L^n_{f - s\tau + zg} \chi_{\alpha})(x_{\alpha}).$$
Now  we repeat the argument in \cite{W} and conclude that
 $$\sum_{|\beta| = m -1} L^n_{f- s\tau + zg} \chi_{\beta} (x_{\beta}) = \sum_{|\alpha| = m}  L^n_{f - s\tau + zg} \chi_{\alpha} (x_{\bar{\alpha}}).$$
Thus the proof of (\ref{eq:2.1}) is reduced to an estimate of the difference
$$L^n_{f - s\tau + zg} \chi_{\alpha}(x_{\alpha}) - L^n_{f - s\tau + zg} \chi_{\alpha}(x_{\bar{\alpha}}).$$
Observe that $x_{\alpha}$ and $x_{\bar{\alpha}}$ are on the same cylinder $U_{\bar{\alpha}}.$ According to Proposition 1, 
the function $L^n_{f - s\tau + zg}\chi_{\alpha}$ is $\nu$-H\"older continuous on $U_{\bar{\alpha}}$. Consequently, for every  $n \geq 2$ we obtain
$$|L^n_{f - s\tau + zg} \chi_{\alpha}(x_{\alpha}) - L^n_{f - s\tau + zg} \chi_{\alpha}(x_{\bar{\alpha}})| 
\leq \|L^n_{f - s\tau + zg}\chi_{\alpha}\|_{\nu} d( x_{\alpha}, x_{\bar{\alpha}})^{\nu},$$
where $\|.\|_{\nu}$ denotes the operator norm derived from the $\nu$-H\"older norm. Going back to (\ref{eq:2.2}), we deduce
\begin{eqnarray} \label{eq:2.3}
\Bigl|Z_n(f - s\tau + zg) - \sum_{i=1}^k L_{f - s\tau + zg}^n \chi_i(x_i)\Bigr| \nonumber \\
\leq \sum_{m=2}^n \sum_{|\alpha| = m} \|L_{f-s\tau +zg}^{n-m}\|_{\nu} \|L^m_{f - s\tau + zg} \chi_{\alpha}\|_{\nu} 
d ( x_{\alpha}, x_{\bar{\alpha}}).
\end{eqnarray}
This it makes possible to apply (\ref{eq:1.1}) and to conclude that
$$d(x_{\alpha}, x_{\bar{\alpha}}) \leq C^{\nu} \gamma_0^{-\nu(m-2)} d(\sigma^{m-2} x_{\alpha}, \sigma^{m-2} 
x_{\bar{\alpha}})^{\nu} \leq C_2 \gamma_0^{-m \nu}.$$
To finish the proof we have to estimate the term $\|L_{g - sr + zf}^m \chi_{\beta}\|_{\nu}.$ Given a word $\alpha$ of length 
$n > 1$ and $x \in \sigma(U_{\alpha_{n-1}}) \cap {\rm Int}\: U_i$, for any $i$ with $A(\alpha_{n-1}, i) = 1$, we define 
$\sigma_{\alpha}^{-1}(x)$ to be the unique point $y$ such that
$\sigma^n(y) = x$ and $y \in U_{\alpha}.$ For a symbol $i$ we define $i x = \sigma_i^{-1} (x).$

First we have
\begin{lem} 
$$(L_{f - s\tau + zg}^m \chi_{\beta})(x) = \begin{cases} e^{(f -s\tau +zg)^m(\sigma_{\beta}^{-1} x)},\: {\rm if}\:x 
\in \sigma(U_{\beta_{m-1}}) ,\\ 0,\:\: {\rm otherwise}. \end{cases}$$
\end{lem}
The proof is a repetition of that of Lemma 3.7 in \cite{W} and it is based on the definition of $\sigma_{\alpha}^{-1}$ above and 
the fact that
$$(L_{f - s\tau + zg}^m\chi_{\beta})(x) = \sum_{\sigma^m y = x} e^{f^m - s\tau^m + zg^m}(y) \chi_{\beta}(y).$$

For every admissible word $\beta$ with $|\beta| = m$, we fix a point $y_{\beta} \in \sigma(U_{\beta_{m-1}})$ which will be chosen 
as in \cite{W}. Define $z_{\beta} = \sigma_{\beta}^{-1} (y_{\beta}).$
\def\sb{\sigma_{\beta}^{-1}}
\begin{lem} There exist constants $B_0 > 0, B_1 > 0, B_2 > 0$ such that we have the estimate
$$\|L_{f - s\tau + zg}^m(\chi_{\beta})\|_{\nu} \leq B_0\Bigl(e^{a_0 |U|^{\nu} B_1} 
+ B_1 |s| e^{a_0 |U|^{\nu} ( 1+ \gamma_0^{-\nu}) B_1}\Bigr)$$
$$\times \Bigl(e^{c_0 |U|^{\nu} B_2} + B_2 |z| e^{c_0 |U|^{\nu} ( 1+ \gamma_0^{-\nu}) B_2}\Bigr)e^{(f^m- a \tau^m 
+ c g^m)(z_{\beta})}.$$
\end{lem}

{\bf Proof.}  We will follow the proof of Lemma 3.8 in \cite{W}. Let  $x$ and $y$ be in  the same Markov leaf. If 
$ y \notin \sigma(U_{\beta_{m-1}})$, then 
$|L_{f - s\tau + zg}^m (\chi_{\beta})(x)| = |L_{f - s\tau + zg}^m (\chi_{\beta})(x)- L_{f- s\tau + zg}^m (\chi_{\beta})(y)| = 0$. 
In the case when $x \notin \sigma(U_{\beta_{m-1}})$, we repeat the same argument. So we will consider the case when both $x$ and 
$y$ are in $\sigma(U_{{\beta}_{m-1}}).$

We have
$$|L_{f - s\tau + zg}^m (\chi_{\beta})(x)| = |e^{(f^m -(a + i b)\tau^m + (c + id) g^m)(\sigma_{\beta}^{-1} x)}|$$
$$\leq \exp\Bigl((f^m - a \tau^m + c g^m)(\sigma_{\beta}^{-1}x) 
-  (f^m - a \tau^m + c g^m)(\sigma_{\beta}^{-1}y)\Bigr)e^{(f^m - a \tau^m + c g^m)(z_{\beta})}.$$
On the other hand, applying Lemma 1 with $w = \tau$, we get
$$|\tau^m(\sigma_{\beta}^{-1} x) - \tau^m(\sigma_{\beta}^{-1}y) | 
\leq B_1 (d(\sigma^{m-1} \sigma_{\beta}^{-1} x, \sigma^{m-1} \sigma_{\beta}^{-1} y))^{\nu} \leq B_1 |U|^{\nu}.$$
The same argument works for the terms involving $f^m$ and $g^m$, applying Lemma 1 with $w= f, g$, respectively. Thus we obtain
$$|L_{f - s\tau + zg}^m (\chi_{\beta})(x)| \leq e^{(C_0 + a_0 B_1 + c_0 B_2)|U|^{\nu}}e^{(f^m - a \tau^m + c g^m)(z_{\beta})}.$$
and this implies an estimate for $|L_{f - s\tau + zg}^m (\chi_{\beta})|_{\infty}.$
Next,
$$|L_{f - s\tau + zg}^m (\chi_{\beta})(x) - L_{f - s\tau + zg}^m(\chi_{\beta})(y)|$$
$$\leq |e^{f^m(\sb(x)) - f^m(\sb(y))} - 1| | e^{f^m (\sb(y))}||e^{-s\tau^m(\sb(x)) + s\tau^m(\sb(y))} - 1 ||e^{-s\tau^m(\sb(y))}|$$
$$\times |e^{zg^m(\sb(x)) - zg^m(\sb(y)) } - 1| |e^{zg^m(\sb(y))}|.$$
As in \cite{W}, we have
$$|e^{-sr^m(\sb(x)) + sr^m(\sb(y))} - 1 ||e^{-sr^m(\sb(y))}| 
\leq B_1 \gamma_0^{\nu} |s|e^{a_0 B_1 (1+ \gamma_0^{-\nu})|U|^{\nu}} e^{-a r^m(z_{\beta})} d(x, y)^{\nu}.$$
For the product involving $zg^m$ we have the same estimate with $B_2, |z|, c_0$ and $c$  in the place of $B_1, |s|, a_0$ and $a$. 
A similar estimate holds for the term containing $f^m$ with a constant $B_3$ in the place of $B_1$. Taking the product of these
estimates, we obtain a bound for $|L_{f - s\tau + zg}^m (\chi_{\beta})(x) - L_{f - s\tau + zg}^m(\chi_{\beta})(y)|$, this implies
the desired estimate for the $\nu$-H\"older norm of $L_{f-m s\tau + zg} (\chi_{\beta})$. This completes the proof. \hfill\qed\\

Now the proof of (\ref{eq:2.1}) is reduced to the estimate of
$$\sum_{|\beta| = m} e^{(f^m - a \tau^m + cg^m)(z_{\beta})}.$$

Introduce the real-valued function $h = f - a\tau + cg$. Then we must estimate 
$$\sum_{|\beta| = m} e^{h^m(z_{\beta})}.$$
For this purpose we repeat the argument on pages 232-234 in \cite{W} and deduce with some constant $d_0 > 0$ depending only on 
the matrix $A$ and every $\epsilon > 0$ the bound
$$\sum_{|\beta| = m} e^{h^m(z_{\beta})} \leq e^{d_0 |h|_{\infty}} B_{\epsilon} e^{(m + d_0) (\epsilon + {\rm Pr}(h))}.$$
Combing this with the previous estimates, we get (\ref{eq:2.1}) and the proof of Theorem 4 is complete. \hfill\qed

\def\tm{\tilde{m}}
\def\tj{\tilde{j}}
\def\dd{{\mathcal D}}
\def\piU{\pi^{(U)}}
\def\hrho{\hat{\rho}}
\def\hdd{\widehat{\dd}}
\def\Xijl{X^{(\ell)}_{i,j}}
\def\hXijl{\widehat{X}^{(\ell)}_{i,j}}
\def\omijl{\omega^{\ell}_{i,j}}
\def\vl{v^{(\ell)}}
\def\kk{{\mathcal K}}
\def\Lip{\mbox{\rm Lip}}
\def\ff{{\mathcal F}}
\def\diamte{\diam_\theta}
\def\J{{\sf J}}
\def\nn{{\mathcal N}}
\def\ma{{\mathcal M}_a}
\def\mac{{\mathcal M}_{ac}}
\def\dte{D_\theta}
\def\fa{f^{(a)}}
\def\lab{{\mathcal L}_{ab}}
\def\dl{d^{(\ell)}}
\def\tu{\tilde{u}}
\def\hu{\hat{u}}
\def\hz{\hat{z}}
\def\Gl{\Gamma^{(\ell)}}
\def\lambdam{\lambda^{(m)}}
\def\thetam{\theta^{(m)}}
\def\tR{\widetilde{R}}
\def\htheta{\hat{\theta}}
\def\hW{\widehat{W}}
\def\hnu{\hat{\nu}}
\def\labw{\lc_{abw}}
\def\labz{\lc_{abz}}
\def\tcc{\tilde{\cc}}
\def\tdd{\widetilde{\dd}}
\def\hcc{\widehat{\cc}}
\def\lambdam{\lambda^{(m)}}
\def\thetam{\theta^{(m)}}
\def\Lipf{\mbox{\rm\footnotesize Lip}}
\def\yl{y^{(\ell)}}
\def\tS{\widetilde{S}}
\def\E{\mathcal E}
\def\tp{\tilde{p}}
\def\tq{\tilde{q}}
\def\f0{f^{(0)}}
\def\ttau{\tilde{\tau}}

\def\fc{f_}
\def\gt{g_t}
\def\tft{\tilde{f}_t}
\def\ft{f_t}
\def\fat{f_{at}}
\def\MM{{\mathcal M}}
\def\mm{{\mathcal M}}
\def\fa{f_a}
\def\fb{f+b}
\def\f0{f_0}
\def\fab{f_{ab}}
\def\mab{\mm_{ab}}
\def\ma{\mm_{a}}
\def\lab{\lc_{ab}}
\def\psib{\psi^{(b)}}
\def\fac{f_{ac}}
\def\mac{\mm_{ac}}
\def\mat{\mm_{at}}
\def\matc{\mm_{atc}}
\def\tmat{\widetilde{\mm}_{at}}
\def\lac{\lc_{ac}}
\def\labc{\lc_{abc}}
\def\labt{\lc_{abt}}
\def\labtz{\lc_{abtz}}
\def\labz{\lc_{abz}}
\def\psic{\psi^{(c)}}

\def\hclip{\hat{C}^{\lip}}
\def\hC{\hat{C}}

\def\PPsi{{\sf \Psi}}
\def\hk{\hat{k}}
\def\hgamma{\hat{\gamma}}

\def\Intu{\mbox{\rm Int}^u}
\def\hZ{\widehat{Z}}
\def\xijl{X^{(\ell)}_{i,j}}
\def\eijl{\omega^{(\ell)}_{i,j}}
\def\hxijl{\widehat{X}^{(\ell)}_{i,j}}

\def\gao{\gamma^{(1)}}
\def\gat{\gamma^{(2)}}
\def\diamtef{{\footnotesize\mbox{\rm  diam}_\theta}}
\def\Con{\mbox{\rm Const}}
\def\gej{\chi^{(j)}_\mu}
\def\ge{\chi_\epsilon}
\def\geo{\chi^{(1)}_\mu}
\def\get{\chi^{(2)}_\mu}
\def\gei{\chi^{(i)}_{\mu}}
\def\gee{\chi_{\mu}}
\def\gett{\chi^{(2)}_{\mu}}
\def\geol{\chi^{(1)}_{\ell}}
\def\getl{\chi^{(2)}_{\ell}}
\def\geil{\chi^{(i)}_{\ell}}
\def\gee{\chi_{\ell}}

\section{Ruelle operators -- definitions and assumptions}
\renewcommand{\theequation}{\arabic{section}.\arabic{equation}}
\setcounter{equation}{0}

\def\chU{\check{U}}

For a contact Anosov flows $\phi_t$ with Lipschitz local stable holonomy maps it is proved in Sect. 6 in
\cite{St2} that the following {\it local non-integrability condition} holds:

\medskip

\noindent
{\sc (LNIC):}  {\it There exist $z_0\in \mt$,  $\ep_0 > 0$ and $\theta_0 > 0$ such that
for any  $\ep \in (0,\ep_0]$, any $\hz\in \mt \cap W^u_{\ep}(z_0)$  and any tangent vector 
$\eta \in E^u(\hz)$ to $\mt$ at $\hz$ with  $\|\eta\| = 1$ there exist  $\tz \in \mt \cap W^u_{\ep}(\hz)$, 
$\ty_1, \ty_2 \in \mt \cap W^s_{\ep}(\tz)$ with $\ty_1 \neq \ty_2$,
$\delta = \delta(\tz,\ty_1, \ty_2) > 0$ and $\ep'= \ep'(\tz,\ty_1,\ty_2)  \in (0,\ep]$ such that
$$|\Delta( \exp^u_{z}(v), \pi_{\ty_1}(z)) -  \Delta( \exp^u_{z}(v), \pi_{\ty_2}(z))| \geq \delta\,  \|v\| $$
for all $z\in W^u_{\ep'}(\tz)\cap\mt$  and  $v\in E^u(z; \ep')$ with  $\exp^u_z(v) \in \mt$ and
$\la \frac{v}{\|v\|} , \eta_z\ra \geq \theta_0$,   where $\eta_z$ is the parallel 
translate of $\eta$ along the geodesic in $W^u_{\ep_0}(z_0)$ from $\hz$ to $z$. }

For any $x \in \mt$, $T > 0$ and  $\delta\in (0,\ep]$ set
$$B^u_T (x,\delta) = \{ y\in W^u_{\ep}(x) : d(\phi_t(x), \phi_t(y)) \leq \delta \: \: , \:\:  0 \leq t \leq T \} .$$

We will say that $\phi_t$ has a {\it regular distortion along unstable manifolds} over
the basic set $\mt$  if there exists a constant $\ep_0 > 0$ with the following properties:

\ms 

(a) For any  $0 < \delta \leq   \ep \leq \ep_0$ there exists a constant $R =  R (\delta , \ep) > 0$ such that 
$$\diam( \mt \cap B^u_T(z ,\ep))   \leq R \, \diam( \mt \cap B^u_T (z , \delta))$$
for any $z \in \mt$ and any $T > 0$.

\ms

(b) For any $\ep \in (0,\ep_0]$ and any $\rho \in (0,1)$ there exists $\delta  \in (0,\ep]$
such that for  any $z\in \mt$ and any $T > 0$ we have
$\diam ( \mt \cap B^u_T(z ,\delta))   \leq \rho \; \diam( \mt \cap B^u_T (z , \ep)) .$

\bs

A large class of flows on basic sets having regular distortion along unstable manifolds  is described in \cite{St3}.

In this paper we work under the following {\bf  Standing Assumptions:} 

\ms

(A) $\phi_t$ has Lipschitz local holonomy maps over $\mt$,

\ms  

(B) the local non-integrability condition (LNIC) holds for $\phi_t$ on $\mt$,

\ms

(C) $\phi_t$  has a regular distortion  along unstable manifolds over the basic set $\mt$.

\ms

A rather large class of examples satisfying the above conditions is provided by imposing the following {\it pinching condition}:

\ms

\noindent
{\bf (P)}:  {\it There exist  constants $C > 0$ and $\beta \geq \alpha > 0$ such that for every  $x\in M$ we have
$$\frac{1}{C} \, e^{\alpha_x \,t}\, \|u\| \leq \| d\phi_{t}(x)\cdot u\| 
\leq C\, e^{\beta_x\,t}\, \|u\| \quad, \quad  u\in E^u(x) \:\:, t > 0 $$
for some constants $\alpha_x, \beta_x > 0$ with
$\alpha \leq \alpha_x \leq \beta_x \leq \beta$ and $2\alpha_x - \beta_x \geq \alpha$ for all $x\in M$.}

\ms

We should note that  (P) holds for geodesic flows on manifolds of strictly negative sectional curvature satisfying the so called $\frac{1}{4}$-pinching condition. 
(P) always holds when  $\dim(M) = 3$.

\ms

{ \bf Simplifying Assumptions:} $\phi_t$ is a $C^2$ contact Anosov flow satisfying the condition (P).

\ms

As shown in \cite{St3}  the pinching condition  (P) implies that $\phi_t$ has  Lipschitz local holonomy maps 
and regular distortion along unstable manifolds. Combining this with Proposition 6.1 in \cite{St2}, shows that the
Simplifying Assumptions imply the Standing Assumptions.
 
 \ms

As in Sect. 1 consider a {\bf fixed Markov family} $\rr = \{ R_i\}_{i=1}^k$ for the flow $\phi_t$ on $\mt$ consisting 
of rectangles $R_i = [U_i,S_i]$ and let  $U = \cup_{i=1}^k U_i$. 
The Standing Assumptions imply  the existence of constants $c_0 \in (0,1]$ and $\gamma_1 > \gamma_0 > 1$ such that
(1.1) hold.

In what follows we will assume that {\bf $f$ and $g$ are fixed  real-valued functions in $C^\alpha(\hU)$  for some fixed $\alpha > 0$}. Let  $P = P_f$ be the unique real 
number so that $\Pr(f- P\, \tau) = 0$, where $\Pr(h)$ is the {\it topological pressure} of $h$ with respect to the shift map $\sigma$ defined in Section 2. 
Given  $t\in \R$ with $t \geq 1$, following \cite{D}, denote by $\ft$ the {\it average of $f$ over balls in $U$ of radius} $1/t$. To be more
precise, first one has to fix an arbitrary extension $f\in C^\alpha(V)$ (with the  same H\"older constant), where $V$ is an open neighborhood of $U$ in $M$, and then take the
averages in question. Then $\ft\in C^\infty(V)$, so its restriction to $U$ is Lipschitz  (with respect to the Riemann metric) and:


(a) $\| f-\ft\|_\infty \leq | f|_\alpha/t^\alpha$ ;
 
(b) $\Lip(\ft) \leq \Con\; \|f\|_\infty t$ ; 

(c) For any $\beta \in (0, \alpha)$ we have
$| f- \ft|_{\beta} \leq 2\, |f|_\alpha/ t^{\alpha-\beta}$.

\ms

In the special case $f\in \clip(U)$ we set $\ft = f$ for all $t \geq 1$. Similarly for $g$. Let $\lambda_0 > 0$ be the largest eigenvalue of 
$L_{f-P\tau}$, and let $\hnu_0$ be the (unique) probability measure on $U$ with $L_{f-P\tau}^*\hnu_0 = \hnu_0$.
Fix a corresponding (positive) eigenfunction $h_{0} \in \hC^\alpha (U)$ such that $\int_U h_{0} \, d\hnu_0 = 1$. Then 
$d\nu_0 = h_0\, d\hnu_0$ defines a {\it $\sigma$-invariant probability measure} $\nu_0$ on $U$. Setting
$$\f0 = f - P\, \tau + \ln h_{0}(u) - \ln h_{0}(\sigma(u)) ,$$
we have $L_{f^{(0)}}^*\nu_0 = \nu_0$, i.e. $\di \int_U L_{f^{(0)}} H \, d\nu_0 = \int_U H\, d\nu_0$ for any $H \in C(U)$, and $L_{\f0} 1 = 1$.

Given real numbers $a$ and $t$ (with $|a| + \frac{1}{|t|}$ small), denote by $\lambda_{at}$ the 
{\it largest eigenvalue} of $L_{\ft -(P+a)\tau}$ on $\clip(U)$ and by $h_{at}$ the corresponding (positive) 
eigenfunction such that $\int_U h_{at}\, d\nu_{at} = 1$, where $\nu_{at}$ is the unique probability measure on  $U$ with $L_{\ft -(P+a)\tau}^*\nu_{at} = \nu_{at}$. 

As is well-known the shift map $\sigma : \hU \longrightarrow \hU$ is naturally isomorphic to an one-sided subshift of finite 
type. Given $\theta \in (0,1)$, a natural metric associated by this isomorphism is defined  (for $x\neq y$) by
$d_\theta(x,y) = \theta^m$, where $m$ is the largest integer such that $x,y$ belong to the same cylinder of length $m$.
There exist $\theta = \theta(\alpha) \in (0,1)$ and $\beta \in (0,\alpha)$ such that $(d(x,y))^\alpha \leq \Con \; d_\theta(x,y)$ 
and $d_\theta(x,y) \leq \Con\; (d(x,y))^\beta$ for all $x,y \in \hU$. 
One can then apply the Ruelle-Perron-Frobenius theorem to the sub-shift of fine type and deduce that $h_{at} \in C^\beta(\hU)$.  However this is not enough 
for our purposes -- in Lemma 4 below we get a bit more.

Consider an arbitrary $\beta \in (0,\alpha)$.
It follows from properties (a) and (c)  above that there exists a constant $C_0 > 0$,
depending on $f$ and $\alpha$ but independent of $\beta$,  such that
\be
\| [\ft - (P+a)\tau] - (f - P \tau)\|_\beta \leq C_0 \, [|a| + 1/t^{\alpha-\beta}]
\ee
for all $|a|\leq 1$ and $t \geq 1$. Since $\Pr(f- P\tau) = 0$, it follows from the analyticity of pressure and the 
eigenfunction projection corresponding to the {\it maximal eigenvalue} $\lambda_{at} = e^{\Prf(\ft - (P+a)\tau)}$ 
of the Ruelle operator $L_{\ft - (P+a)\tau}$ on $C^\beta(U)$ (cf. e.g. Ch. 3 in \cite{PP}) that there exists a constant
$a_0 > 0$  such that, taking $C_0 > 0$ sufficiently large, we have
\be
|\Pr(\ft - (P+a)\tau)| \leq C_0 \; \left(|a| + \frac{1}{t^{\alpha-\beta}}\right) \quad , \quad
\|h_{at} - h_0\|_\beta \leq  C_0 \; \left(|a| + \frac{1}{t^{\alpha-\beta}}\right)\;
\ee
for  $|a| \leq a_0$ and $1/t \leq a_0$. We may assume $C_0 > 0$  and $a_0> 0$ are taken so that $1/C_0 \leq \lambda_{at} \leq C_0$,
$\|\ft\|_\infty \leq C_0$ and $1/C_0 \leq h_{at}(u) \leq C_0$ for all $u \in U$ and all $|a|, 1/t \leq a_0$.

Given real numbers $a$ and $t$ with $|a|, 1/t \leq a_0$ consider the functions
$$\fat = \ft - (P+a) \tau + \ln h_{at} - \ln (h_{at}\circ \sigma) - \ln \lambda_{at}$$
and the operators 
$$ \labt = L_{\fat - \i\, b\, \tau} : C(U) \longrightarrow C(U)\:\:\: , \:\:\:
\mat = L_{\fat} : C(U) \longrightarrow C(U).$$
One checks that $\mat \; 1 = 1$.

Taking the constant $C_0 > 0$ sufficiently large, we may assume that
\be
\|\fat - \f0\|_\beta \leq C_0\, \left[ |a| +  \frac{1}{t^{\alpha- \beta}} \right] \quad, \quad
|a|, 1/t \leq a_0 .
\ee

We will now prove a simple uniform estimate for $\Lip(h_{at})$. With respect to the usual metrics on
symbol spaces this a consequence of general facts (see e.g. Sect. 1.7 in \cite{B1} or Ch. 3 in \cite{PP}),
however here we need it with respect to the Riemann metric.

 The proof of the following lemma is given in the Appendix.

\begin{lem} Taking the constant $a_0 > 0$ sufficiently small, there exists a constant $T' > 0$ such that 
for all $a,t \in \R$ with $|a|\leq a_0$ and $t \geq 1/a_0$ we have $h_{at} \in \clip(\hU)$ and
$\Lip(h_{at}) \leq T' t$.
\end{lem}

\ms

It follows from the above that, assuming $a_0 > 0$ is chosen sufficiently small, there exists a constant $T > 0$ 
(depending on $|f|_\alpha$ and $a_0$) such that 
\be
\|\fat\|_\infty \leq T \quad , \quad \|\gt\|_\infty \leq T \quad , \quad\Lip(h_{at}) \leq T\, t
\quad , \quad \Lip(\fat) \leq T\, t
\ee
for $|a|, 1/t \leq a_0$. We will also assume that $T \geq \max \{ \, \|\tau \|_0 \, , \, \Lip(\tau_{|\hU}) \, \}.$
From now on {\bf we will assume that $a_0$, $C_0$, $T$, $1 < \gamma_0 < \gamma_1$ are fixed constants}  
with (1.1) and (3.1) -- (3.4).

\section{Ruelle operators depending on two parameters -- the case when $b$ is the leading parameter}
\renewcommand{\theequation}{\arabic{section}.\arabic{equation}}
\setcounter{equation}{0}

Throughout this section we work under the Standing Assumptions made in Sect. 3 and with fixed
real-valued functions $f,g \in C^\alpha(\hU)$ as in Sect. 3. Throughout $0 < \beta < \alpha$ are fixed numbers.

We will study Ruelle operators of the form $L_{f - (P_f + a + \i b)\tau + z g}$, where $z = c + \i w$, 
$a,b, c, w\in \R$, and $|a|, |c|\leq a_0$  for some constant $a_0 > 0$. Such
operators will be approximates by operators of the form
$$\labtz = L_{\fat - \i\,b\tau + z \gt} : C^{\alpha}(\hU) \longrightarrow C^\alpha (\hU) .$$
In fact, since $\fat - \i b \tau + z \gt$ is Lipschitz, the operators $\labtz$ preserves
each of the spaces $C^{\alpha'}(\hU)$ for $0 < \alpha' \leq 1$ including the space
$\clip (\hU)$ of Lipschitz functions $h: \hU \longrightarrow \C$. For such $h$ we will denote by 
$\Lip(h)$ the Lipschitz constant of $h$. Let $\| h\|_0$ denote the {\it standard $\sup$ norm}  
of $h$ on $\hU$.  For  $|b| \geq 1$, as in \cite{D}, consider the norm $\|.\|_{\lip,b}$ on $\clip (\hU)$ defined by 
$\| h\|_{\lip,b} = \|h\|_0 + \frac{\lip(h)}{|b|}$. and also the norm
$\|h\|_{\beta,b} = \|h\|_\infty + \frac{|h|_\beta}{|b|}$
on $C^\beta(U)$.

Our aim in this  section is to prove the following

\begin{thm} 
Let $\phi_t : M \longrightarrow M$ satisfy the Standing Assumptions over the basic set $\mt$, and let $0 < \beta < \alpha$.
Let  $\rr = \{R_i\}_{i=1}^k$ be a Markov family for $\phi_t$ over $\mt$ as in Sect. 1. Then 
for any real-valued functions  $f,g \in C^\alpha(\hU)$ we have:

\ms

{\rm (a)} {\it For any constants $\epsilon > 0$, $B > 0$ and $\nu \in (0,1)$ there exist constants $0 < \rho < 1$, $a_0 > 0$, 
$b_0 \geq 1$, $A_0 > 0$ and  $C = C(B, \epsilon)> 0$ such that if $a,c\in \R$  satisfy $|a|, |c| \leq a_0$,   then
$$\|L_{\fat - \i b \tau + (c+\i w) \gt}^m h \|_{\lip,b}  \leq C \, \rho^m \, |b|^{\ep}\, \| h\|_{\lip,b}$$
for all $h\in \clip(\hU)$, all  integers $m \geq 1$ and all $b, w, t\in \R$ with  $|b| \geq b_0$,
$1 \leq t \leq \frac{1}{A_0} \log|b|^\nu$ and $|w| \leq B \, |b|^\nu$.}

\ms

{\rm (b)} {\it For any constants $\epsilon > 0$, $B > 0$, $\nu \in (0,1)$ and $\beta\in (0,\alpha)$  
there exist constants $0 < \rho < 1$, $a_0 > 0$, $b_0 \geq 1$ and  $C = C(B, \epsilon)> 0$ such that if $a,c\in \R$  
satisfy $|a|, |c| \leq a_0$, then
$$\|L_{f -(P_f+a+ \i b)\tau + (c+\i w) g}^m h \|_{\beta,b}  \leq C \, \rho^m \, |b|^{\ep}\, \| h\|_{\beta,b}$$
for all $h \in C^\beta(\hU)$, all integers $m \geq 1$ and all $b, w\in \R$ with  $|b| \geq b_0$ and $|w| \leq B \, |b|^\nu$.}

\ms

{\rm (c)} {\it If $f,g \in \clip(\hU)$, then for any constants $\epsilon > 0$, $B > 0$ and $\beta\in (0,\alpha)$  
there exist constants $0 < \rho < 1$, $a_0 > 0$, $b_0 \geq 1$ and  $C = C(B, \epsilon)> 0$ such that if $a,c\in \R$  
satisfy $|a|, |c| \leq a_0$, then
$$\|L_{f -(P_f+a+ \i b)\tau + (c+\i w) g}^m h \|_{\lip,b}  \leq C \, \rho^m \, |b|^{\ep}\, \| h\|_{\lip,b}$$
for all $h \in C^\beta(\hU)$, all integers $m \geq 1$ and all $b, w\in \R$ with  $|b| \geq b_0$ and $|w| \leq B \, |b|$.}

\end{thm}

\medskip
We will first prove part (a) of the above theorem and then derive part (b) by a simple approximation procedure.
To prove part (a) we will use the main steps in Section 5 in \cite{St2} with necessary modifications. 
The proof of part (c) is just a much simpler version of the proof of (b).

Define a {\it new metric} $D$ on $\hU$ by 
$$D(x,y) = \min \{ \diam(\cc) : x,y\in \cc\:, \: \cc \: \mbox{\rm a cylinder contained in }\, U_i \}$$
if $x,y \in U_i$ for some $i = 1, \ldots,k$, and $D(x,y) = 1$ otherwise. Rescaling the metric on $M$ if necessary,
we will assume that $\diam(U_i) < 1$ for all $i$. As shown in \cite{St1}, $D$ is a metric on $\hU$ with $d(x,y) \leq  D(x,y)$ 
for $x,y\in \hU_i$ for some $i$, and for any cylinder $\cc$ in $U$ the characteristic function  $\chi_{\hcc}$ of $\hcc$ on 
$\hU$ is Lipschitz with respect to $D$ and $ \Lip_D(\chi_{\hcc}) \leq 1/\diam(\cc)$.

We will denote by $\clip_D(\hU)$ the {\it space of all Lipschitz functions $h : \hU \longrightarrow \C$
with respect to the metric} $D$ on $\hU$ and by $\Lip_D(h)$ the {\it Lipschitz constant} of $h$ with respect to $D$.

Given $A > 0$, denote by $K_A(\hU)$  {\it the set of all functions } $h\in \clip_D(\hU)$  such that  $h > 0$ and
$\frac{|h(u) - h(u')|}{h(u')} \leq A\, D (u,u')$ for all $u,u' \in \hU$ that belong to the  same $\hU_i$  for some $i = 1, \ldots,k$.
Notice that $h\in K_A(\hU)$ implies $|\ln h(u) - \ln h(v)| \leq A\; D (u,v)$ and therefore
$e^{-A\; D (u,v)} \leq \frac{h(u)}{h(v)} \leq e^{A \; D (u,v)}$ for all $u, v\in \hU_i $, $i = 1, \ldots,k.$

We begin with a lemma of Lasota-Yorke type, which necessarily has a more complicated form due to the
more complex situation considered. It involves the operators $\labtz$, and also operators of the form
$$  \matc = L_{\fat + c \gt} : C^\alpha (\hU) \longrightarrow C^\alpha (\hU) .$$
{\bf Fix arbitrary constants $\nu \in (0,1)$ and $\hgamma$ with} $ 1 < \hgamma < \gamma_0$.

\ms

\begin{lem} Assuming $a_0 > 0$ is chosen sufficiently small, there exists a constant $A_0 > 0$ 
such that for all $a,c, t\in \R$ with $|a|, |c|\leq a_0$ and $t \geq 1$ the following hold:

\ms

(a)  {\it If $H \in K_E(\hU)$ for some $E > 0$, then 
$$\frac{|(\matc^m H)(u) - (\matc^m H)(u')|}{(\matc^m H)(u')} \leq
A_0 \, \left[ \frac{E}{\hgamma^m} + e^{A_0 t} \, t \right]\, D (u,u')$$
for all $m \geq 1$ and all $u,u'\in U_i$, $i = 1, \ldots, k$.}

\ms

(b) {\it If the functions $h$ and  $H$ on $\hU$  and $E > 0$  are such that $H > 0$ on $\hU$ and 
$|h(v) - h(v')| \leq E\, H(v')\, D (v,v')$ for any $v,v'\in \hU_i$, $i = 1, \ldots,k$, 
then for any integer $m \geq 1$ and any $b, w, t\in \R$ with  $|b|, t, |w|\geq 1$, for $z = c + \i w$  we have 
$$ |\labtz^m h(u) - \labtz^m h(u')| \leq  
 A_0  \left(\frac{E}{\hgamma^m} (\matc^m H)(u') + ( |b| + e^{A_0t}t + t|w|) (\matc^m |h|)(u')\right)\, D(u,u') $$
whenever $u,u'\in \hU_i$ for some $i = 1, \ldots,k$. In particular, if 
\be
t \leq \frac{\log |b|^\nu}{A_0} \quad , \quad t \leq B |b|^{1-\nu}\quad , \quad \quad |w| \leq B |b|^\nu
\ee
for some constant $B > 0$, then 
$$ |\labtz^m h(u) - \labtz^m h(u')| \leq  
 A_1  \left(\frac{E}{\hgamma^m} (\matc^m H)(u') + |b| (\matc^m |h|)(u')\right)\, D(u,u') .$$
 for some constant $A_1 > 0$. }
\end{lem}

A proof of this lemma is given in the Appendix.

{\bf From now on we will assume that $a_0$, $\eta_0$ and $A_0$ are fixed with the properties in Lemma 5 above and
$a,b,c,w,t\in \R$ are such that $|a| \leq a_0$, $c \leq \eta_0$, $|b|, t, |w| \geq 1$ and (4.1) hold. As before, 
set $z = c+ \i d$.}

We will use the entire set-up and notation from Section 5 in \cite{St2}. In what follows we recall  the main part of it.

Following Sect. 4 in \cite{St2}, {\bf fix  an arbitrary point $z_0 \in \mt$ and constants  $\ep_0 > 0$ and $\theta_0 \in (0,1)$ 
 with the properties described  in} (LNIC). Assume that  $z_0 \in \Int_\mt(U_1)$, $U_1 \subset \mt \cap W^u_{\ep_0}(z_0)$ and 
$S_1 \subset \mt \cap W^s_{\ep_0}(z_0)$. Fix an arbitrary constant $\theta_1$ such that
$$0 <  \theta_0  < \theta_1 < 1 \;. $$

Next, fix an arbitrary orthonormal basis $e_1, \ldots, e_{n}$ in $E^u (z_0)$ and a $C^1$
parametrization $r(s) = \exp^u_{z_0}(s)$, $s\in V'_0$, of a small neighborhood $W_0$ of $z_0$ in 
$W^u_{\ep_0} (z_0)$ such that $V'_0$ is a convex compact neighborhood of $0$ in 
$\R^{n} \approx \mbox{\rm span}(e_1, \ldots,e_n) = E^u(z_0)$. Then $r(0) = z_0$ and
$\frac{\partial}{\partial s_i} r(s)_{|s=0} = e_i$ for all $i = 1, \ldots,n$.  Set  $U'_0 = W_0\cap \mt $.
Shrinking $W_0$ (and therefore $V'_0$ as well)
if necessary, we may assume that $\overline{U'_0} \subset \Int_\mt (U_1)$ and
$\left|\left\la \frac{\partial r}{\partial s_i} (s) ,  \frac{\partial r}{\partial s_j} (s) \right\ra
- \delta_{ij}\right| $ is uniformly small for all $i, j = 1, \ldots, n$ and $s\in V'_0$, so that
$$\frac{1}{2} \la \xi , \eta \ra \leq \la \; d r(s)\cdot \xi \; , \;  d r(s)\cdot \eta\;\ra \leq 2\, \la \xi ,\eta \ra \quad, \quad
\xi , \eta \in E^u(z_0) \:,\: s\in V'_0\, ,$$
and
$\frac{1}{2}\, \|s-s'\| \leq  d ( r(s) ,  r(s')) \leq 2\, \|s-s'\|$, $s,s'\in V'_0$.

\ms

\noindent
{\bf Definitions} (\cite{St2}): (a)
For a cylinder $\cc \subset U'_0$ and a unit vector $\xi \in E^u(z_0)$
we will say that a {\it separation by a $\xi$-plane occurs} in $\cc$ if there exist $u,v\in \cc$ with 
$d(u,v) \geq \frac{1}{2}\, \diam(\cc)$ such that
$ \left\la \frac{r^{-1}(v) - r^{-1}(u)}{\| r^{-1}(v) - r^{-1}(u)\|}\;,\; \xi \right\ra  \geq \theta_1\;.$

Let  $\ss_\xi$ be the {\it family of all cylinders} $\cc$ contained in $U'_0$ such that a separation by an $\xi$-plane 
occurs in $\cc$.

\ms

(b) Given an open subset $V$ of $U'_0$  which is a finite union of open cylinders and  $\delta > 0$, let
$\cc_1, \ldots, \cc_p$ ($p = p(\delta)\geq 1$) be the family of maximal closed cylinders in $\oV$ with
$\diam(\cc_j) \leq \delta$. For any unit vector $\xi \in E^u(z_0)$ set
$M_{\xi}^{(\delta)}(V) = \cup \{ \cc_j : \cc_j \in \ss_{\xi} \:, \: 1\leq j \leq p\}\;.$

\ms

In what follows we will construct, amongst other things, a sequence of unit vectors
$\xi_1, \xi_2, \ldots, \xi_{j_0}\in E^u(z_0)$. For each $\ell = 1, \ldots,j_0$  set 
$B_\ell = \{ \eta \in \sn : \la \eta , \xi_\ell\ra  \geq \theta_0\}\;.$ For $t \in \R$ and $s\in E^u(z_0)$ set
$I_{\eta ,t} g(s) = \frac{g(s+t\, \eta) - g(s)}{t}$, $ t \neq 0\;$ ({\it increment} of $g$ in the direction of $\eta$).

\begin{lem}
{\rm (\cite{St2})}
There exist integers $1 \leq n_1 \leq N_0$ and $\ell_0 \geq 1$,
a sequence of unit vectors $\eta_1, \eta_2, \ldots, \eta_{\ell_0}\in E^u(z_0)$
and a non-empty open subset $U_0$ of $U'_0$ which is a finite union of open cylinders of 
length $n_1$ such that setting $\uu = \sigma^{n_1} (U_0)$ we have:

{\rm (a)} {\it For any integer $N\geq N_0$ there exist Lipschitz maps $\vl_1, \vl_2 : U \longrightarrow U$ 
($\ell = 1,\ldots, \ell_0$)  such that $\sigma^N(\vl_i(x)) = x$  for all $x\in \uu$ and $\vl_i (\uu)$ is
a finite union of open cylinders of length $N$ ($i=1,2$; $\ell = 1,2, \ldots,\ell_0$).} 

{\rm (b)} {\it There exists a constant $\hd > 0$ such that for all  $\ell = 1, \ldots, \ell_0$, 
$s\in r^{-1}(U_0)$, $0 < |h| \leq \hd$ and $\eta \in B_\ell$ with  $s+h\, \eta \in r^{-1}(U_0\cap \mt)$ we have }
$$\left[I_{\eta,h} \left(\tau^{N}(\vl_2(\trr(\cdot ))) - \tau^{N}(\vl_1(\trr(\cdot)))\right)\right](s)  \geq \frac{\hd}{2}\,.$$

{\rm (c)} {\it We have $\overline{\vl_i (U)} \bigcap \overline{v_{i'}^{(\ell')}(U)} = \e$ whenever $(i,\ell) \neq (i',\ell')$.}

{\rm (d)} {\it For  any open cylinder $V$ in  $U_0$ there exists a constant  $\delta' = \delta'(V) > 0$  such that
$$V \subset M_{\eta_1}^{(\delta)}(V) \cup M_{\eta_2}^{(\delta)}(V) \cup \ldots  \cup M_{\eta_{\ell_0}}^{(\delta)}(V)$$ 
for all} $\delta \in (0,\delta'].$
\end{lem}


Fix  $U_0$ and $\uu$ with the properties described in Lemma 1; then $\overline{\uu} = U$.

\def\Ulo{U^{(\ell_0)}}

Set  $\di \hd = \min_{1\leq \ell\leq \ell_0} \hd_j$, $\di n_0 = \max_{1\leq \ell\leq \ell_0} m_\ell$, 
and fix an arbitrary point $\hz_0 \in U_0^{(\ell_0)}\cap \hU$. 

Fix integers $1 \leq n_1 \leq N_0$ and $\ell_0 \geq 1$, unit vectors $\eta_1, \eta_2, \ldots, \eta_{\ell_0}\in E^u(z_0)$
and a non-empty open subset $U_0$ of $W_0$ 
with the properties described in Lemma 6.
By the choice of $U_0$, $\sigma^{n_1} : U_0 \longrightarrow \uu$ is one-to-one and has an inverse
map $\psi : \uu \longrightarrow U_0$, which is Lipschitz.


Set $E = \max\left\{  4A_0\;, \; \frac{2 A_0\, T}{\gamma-1}\; \right\},$
where $A_0 \geq 1$ is the constant from Lemma 5.4, and {\bf fix an integer} $N \geq N_0$ such that
$$\gamma^N \geq \max \left\{ \: 6A_0 \; , \; \frac{200\, \gamma_1^{n_1}\,A_0}{c_0^2} \; , \; 
\frac{512\, \gamma^{n_1}\, E}{c_0\, \hd\, \rho}  \; \right\} \, .$$
Then fix maps $\vl_i : U \longrightarrow U$ ($\ell = 1, \ldots, \ell_0$, $i = 1,2$) with
the properties (a), (b), (c) and (d) in Lemma 6.  In particular, (c) gives
$$\overline{v_{i}^{(\ell)}(U)} \cap \overline{v_{i'}^{(\ell')}(U)} = \e \quad ,\quad (i,\ell) \neq (i',\ell') .$$

Since $U_0$ is a finite union of open cylinders, it follows from Lemma 6(d)  that there exist a  constant
$\delta' = \delta'(U_0) > 0$  such that
$$M_{\eta_1}^{(\delta)}(U_0) \cup \ldots\cup M_{\eta_{\ell_0}}^{(\delta)}(U_0) \supset U_0  \quad , \quad \delta\in (0, \delta'] .$$
{\bf Fix $\delta'$ with this property}. Set
$$\ep_1 = \min\left\{ \;\frac{1}{32 C_0 }\;,\; c_1\;,\;  \frac{1}{4E}\;,\;
\frac{1}{\hd\, \rho^{p_0+2} } \; , \; \frac{c_0r_0}{\gamma_1^{n_1}}\; , \; \frac{c_0^2(\gamma-1)}{16T\gamma_1^{n_1}}\, \right\} ,$$
and let $b\in \R$ be such that $|b| \geq 1$ and
 $$\frac{\ep_1}{|b|} \leq \delta' .$$

Let $\cc_m$  ($1\leq m \leq p$) be the family of {\it maximal  closed cylinders} 
contained in $\overline{U_0}$ with  $\diam(\cc_m )\leq \frac{\ep_1}{|b|}$ such that 
$U_0 \subset \cup_{j=m}^p \cc_m$ and $\overline{U_0} = \cup_{m=1}^p \cc_m$. As in \cite{St2},
\be
 \rho\, \frac{\ep_1}{|b|} \leq \diam(\cc_m ) \leq  \frac{\ep_1}{|b|} \quad, \quad 1\leq m\leq p\;.
\ee

Fix an integer $q_0 \geq 1$ such that
$$\theta_0 < \theta_1  - 32\, \rho^{q_0-1} .$$
Next, let $\dd_1, \ldots, \dd_q$  be the list of all closed cylinders contained in $\overline{U_0}$ 
that are {\it subcylinders of co-length} $p_0\, q_0$ of some $\cc_m$ ($1\leq m\leq p$). 
Then  $\overline{U_0} = \cc_1 \cup \ldots \cup \cc_p =  \dd_1 \cup \ldots \cup \dd_q$.
Moreover,
$$\rho^{p_0\, q_0+1}\cdot \frac{\ep_1}{|b|} \leq \diam(\dd_j ) \leq  \rho^{q_0}\cdot  \frac{\ep_1}{|b|} \quad, \quad 1\leq j\leq q.$$

Given $j = 1, \ldots,q $, $\ell = 1, \ldots, \ell_0$  and $i = 1,2$, set $\hdd_j = \dd_j \cap \hU$, 
$Z_j = \overline{\sigma^{n_1}(\hdd_j)}$, $\hZ_j = Z_j \cap \hU$, $\xijl = \overline{\vl_i(\hZ_j)}$,
and $\hxijl = \xijl \cap \hU$.  It then follows that $\dd_j =  \psi(Z_j)$, 
and $U = \cup_{j=1}^q Z_j$.  Moreover, $\sigma^{N-n_1}(\vl_i(x)) = \psi(x)$ for all $x \in \uu$, and
all $\xijl$ are cylinders such that $\xijl \cap {X}_{i', j'}^{(\ell')} = \e$  whenever $(i,j,\ell) \neq  (i',j', \ell')$, and
$$\diam(\xijl) \geq \frac{c_0\, \rho^{p_0\, q_0+1}}{\gamma_1^{N}} \cdot \frac{\ep_1}{|b|}$$
for all $i = 1,2$, $j = 1, \ldots,q $ and $\ell = 1, \ldots, \ell_0$. The {\it characteristic function} 
$\omega_{i,j}^{(\ell)} = \chi_{\hxijl} : \hU \longrightarrow [0,1]$ of $\hxijl$
belongs to $\clip_D(\hU)$ and $\Lip_D(\xijl) \leq 1/\diam(\xijl)$. 

Let $J$  be a {\it subset of the set} 
$\Xi = \{\; (i,j, \ell) \; :  \; 1\leq i \leq 2\; ,\;  1\leq j\leq q\; , \;  1\leq \ell \leq \ell_0\;\}.$
Set
$$\mu_0 = \mu_0 (N) = \min \left\{\; \frac{1}{4} \; , \;  \frac{c_0  \, \rho^{p_0q_0+2}\, \ep_1}{4\, \gamma_1^N }\; , \;
\frac{1}{4\,e^{2 T N}} \,\sin^2\left(\frac{\hd\, \rho\, \ep_1}{256}\right) \; \right\} ,$$
and define the function  $\omega = \omega_{J} : \hU \longrightarrow [0,1]$ by
$\di \omega  = 1- \mu_0 \,\sum_{(i, j, \ell) \in J} \eijl\;.$
Clearly $\omega \in \clip_D(\hU)$ and $1-\mu \leq \omega(u) \leq 1$ for any $u \in \hU$. Moreover,
$$\Lip_{D} (\omega) \leq \Gamma =  \frac{ 2 \mu \,\gamma_1^N}{c_0\,\rho^{p_0q_0+2}}\cdot \frac{|b|}{\ep_1} .$$

Next, define {\it the contraction operator} $\nn = \nn_J(a,b,t,c) : \clip_D (\hU) \longrightarrow \clip_D (\hU)$ by 
$$\left(\nn h\right) = \matc^N (\omega_J \cdot h) .$$

Using Lemma 5 above, the proof of the following lemma is the same as that of Lemma 5.6 in \cite{St2}.

\begin{lem} 
{\it Under the above conditions for $N$ and $\mu$ the following hold :}

(a) $\nn h\in K_{E|b|}(\hU)$ {\it for any} $h\in K_{E|b|}(\hU)$;

(b) {\it If $h \in \clip_D (\hU)$ and $H \in K_{E|b|}(\hU)$ are such that $|h|\leq H$ in $\hU$ and
$| h (v) - h(v')|\leq E|b| H(v')\, D (v,v')$
for any $v,v'\in U_j$, $j = 1,\ldots,k$, then for any $i = 1, \ldots,k$ and any $u,u'\in \hU_i $ we have}
$$| (\labtz^N  h)(u) - (\labtz^N h)(u')| \leq E |b| (\nn H)(u')\, D (u,u') .$$
\end{lem}

\ms

\noindent
{\bf Definition.} 
A subset $J$ of $\Xi$ will be called {\it dense} if for any $m = 1,\ldots, p$ 
there exists $(i,j, \ell)\in J$ such that $\dd_j \subset \cc_m$.

\ms

Denote by $\J = \J(a,b)$ the {\it set of all dense subsets} $J$ of $\Xi$.

Although the operator $\nn$ here is different, the proof of the following lemma is very similar to that of 
Lemma 5.8 in \cite{St2}.

\begin{lem}
Given the number $N$, there exist  $\rho_2 = \rho_2(N) \in (0,1)$ and $a_0 = a_0(N) > 0$ 
such that $\di\int_{\hU} (\nn_J H)^2d\nu \leq \rho_{2} \,\int_{\hU} H^2 d\nu$ whenever 
$|a|, |c| \leq a_0$, $t \geq 1/a_0$, $J$ is dense and $H \in K_{E|b|}(\hU)$.
\end{lem}

In what follows we assume that $h, H\in \clip_D (\hU)$ are such that 
\begin{equation}
H\in K_{E|b|}(\hU)\quad , \quad |h(u)|\leq H(u) \:\:\:\:, \:\:\; u\in \hU \;,
\ee
and
\be
|h(u) - h(u')|\leq E|b| H(u')\, D (u,u')\:\:\:\:\:
\mbox{\rm whenever} \:\: u,u'\in \hU_i\;, \; i = 1, \ldots,k \;.
\end{equation}  

Let again $z = c + \i w$.
Define the functions $\geil : \hU  \longrightarrow \C$ ($\ell = 1, \ldots, j_0$, $i = 1,2$) by
$$\di \geol(u) = \frac{\di \left| e^{(\fat^{N} - \i b\tau^{N} + z \gt^N)(\vl_1(u))} h(\vl_1(u)) +
 e^{(\fat^{N} - \i b\tau^{N} + z \gt^N)(\vl_2(u))} h(\vl_2(u))\right|}{\di (1-\mu)
e^{\fat^{N}(\vl_1(u)) + c \gt^N( \vl_1(u))}H(\vl_1(u)) + e^{\fat^{N}(\vl_2(u))+ c \gt^N(\vl_2(u))}H(\vl_2(u))} ,$$ 
$$\di \getl(u) = \frac{\di \left| e^{(\fat^{N} - \i b\tau^{N} + z \gt^N)(\vl_1(u))} h(\vl_1(u)) +
 e^{(\fat^{N} - \i b\tau^{N} + z \gt^N)(\vl_2(u))} h(\vl_2(u))\right|}{\di
e^{\fat^{N}(\vl_1(u)) + c \gt^N(\vl_1(u))}H(\vl_1(u)) +  (1-\mu) e^{\fat^{N}(\vl_2(u)) + c\gt^N(\vl_2(u))}H(\vl_2(u))} ,$$ 
and set $ \gl(u) = b\, [\tau^N(\vl_2(u)) - \tau^N(\vl_1(u))]$, $u\in \hU$.
 
\bs

\noindent
{\bf Definitions.} 
We will say that the cylinders $\dd_j$ and $\dd_{j'}$ are {\it adjacent} if they 
are subcylinders of the same $\cc_m$ for some $m$. 
If $\dd_j$ and $\dd_{j'}$ are contained in $\cc_m$  for some $m$
and  for some  $\ell = 1, \ldots, \ell_0$ there exist  $u \in \dd_j$ and $v\in \dd_{j'}$  
such that $d(u,v) \geq \frac{1}{2}\, \diam(\cc_m)$ and 
$\left\la \frac{r^{-1}(v) - r^{-1}(u)}{\| r^{-1}(v) - r^{-1}(u)\|}\;,\; \eta_\ell \right\ra  \geq \theta_1$,
we will say that $\dd_j$ and $\dd_{j'}$ are {\it $\eta_\ell$-separable in $\cc_m$}.

As a consequence of Lemma 6(b) one gets the following.

\begin{lem}
{\rm (Lemma 5.9 in \cite{St2})}
Let $j, j'\in \{ 1, 2,\ldots,q\}$ be
such that $\dd_j$ and $\dd_{j'}$ are contained in $\cc_m$ and  are $\eta_\ell$-separable 
in $\cc_m$ for some $m = 1, \ldots, p$ and $\ell = 1, \ldots, \ell_0$ . Then
$ |\gl (u) - \gl(u')| \geq c_2\epsilon_{1}$ for all $u\in \hZ_j$ and  $u'\in \hZ_{j'}$,
where $\di c_2 =  \frac{\hd\, \rho}{16}$.
\end{lem}

The  following lemma is the analogue of Lemma 5.10 in \cite{St2} and represents the main step in proving Theorem 1.

\begin{lem}
Assume $|b| \geq b_0$ for some sufficiently large $b_0 > 0$, $|a|, |c| \leq a_0$, and let {\rm (4.1)} hold. 
Then  for any $j = 1, \ldots,q$  there exist $i \in \{ 1,2\}$, $j' \in \{ 1,\ldots,q\}$ and 
$\ell \in \{ 1, \ldots, \ell_0\}$  such that  $\dd_j$ and $\dd_{j'}$ 
are adjacent and $\chi^{(i)}_{\ell} (u) \leq 1$ for all  $u\in \hZ_{j'}$. 
\end{lem}

To prove this  we need the following lemma which coincides with Lemma 14 in \cite{D} and its proof is almost the same.

\begin{lem}
 If $h$ and $H$ satisfy {\rm (4.3)-(4.4)}, then for any $j = 1, \ldots,q$, $i = 1,2$ and $\ell = 1,\ldots,  \ell_0$ we have:

(a) {\it $\di\frac{1}{2} \leq \frac{H(\vl_i(u'))}{H(\vl_i(u''))} \leq 2$ for all} $u', u'' \in \hZ_j$;

(b) {\it Either for all $u\in \hZ_j$ we have $|h(\vl_i(u))|\leq \frac{3}{4}H(\vl_i(u))$, or  
$|h(\vl_i(u))|\geq \frac{1}{4}H(\vl_i(u))$
for all $u\in \hZ_j$.}
\end{lem}

\ms

\noindent
{\it Sketch of proof of Lemma } 10. We use a modification of the proof of Lemma 5.10 in \cite{St2}.

Given $j = 1, \ldots, q$, let $m = 1, \ldots, p$ be such that $\dd_j \subset \cc_m$.  
As in \cite{St2} we find $j',j'' = 1, \ldots,q$ such that $\dd_{j'}, \dd_{j''} \subset \cc_m$ and 
$\dd_{j'}$ and $\dd_{j''}$ are $\eta_\ell$-separable in $\cc_m$. 

Fix $\ell$, $j'$ and $j''$ with the above properties, and set $\hZ = \hZ_j \cup \hZ_{j'}\cup \hZ_{j''}\; .$ 
If there exist $t \in \{j, j', j''\}$ and $i = 1,2$ such 
that the first alternative in Lemma 11(b) holds for $\hZ_{t}$, $\ell$  and $i$, then $\mu \leq 1/4$
implies $\chi_\ell^{(i)}(u) \leq 1$ for any $u\in \hZ_{t}$.

Assume that for every $t\in \{ j, j', j''\}$ and every $i = 1,2$ the second alternative 
in Lemma 11(b) holds for $\hZ_{t}$, $\ell$ and $i$, i.e.  $|h(\vl_i(u))|\geq \frac{1}{4}\, H(\vl_i(u))$, $u \in \hZ$.

Since $\psi(\hZ) = \hdd_j \cup \hdd_{j'} \cup \hdd_{j''} \subset \cc_m$, given $u,u'\in \hZ$ we have 
$\sigma^{N-n_1}(\vl_i(u)), \sigma^{N-n_1}(\vl_i(u'))\in \cc_m$.  Moreover, $\cc' = \vl_i (\sigma^{n_1}(\cc_m))$ is a cylinder with
$\diam(\cc') \leq \frac{\ep_1}{c_0\, \gamma^{N-n_1}\, |b|}$. 
Thus, 
the estimate (8.3) in the Appendix below implies
$$ |\gt^N(\vl_i(u)) - \gt^N(\vl_i(u'))| \leq  \frac{C_1 t \ep_1}{c_0\, \gamma^{N-n_1}\, |b|} .$$
Using the above assumption, (4.1), (4.2) and (3.5),
and assuming e.g.   
$$e^{c \gt^N(\vl_i(u))} |h(\vl_i(u))| \geq e^{c \gt^N(\vl_i(u))} |h(\vl_i(u'))| ,$$ 
we get\footnote{Using some estimates as in the proof of Lemma 5(b) in the Appendix below and 
$\|c \gt^N\|_0 \leq a_0 NT$ by (3.5).}
\begin{eqnarray*}
&      & \frac{|e^{z \gt^N(\vl_i(u))} h(\vl_i(u)) 
- e^{z \gt^N(\vl_i(u'))} h(\vl_i(u'))|}{\min\{| e^{z \gt^N(\vl_i(u))} h(\vl_i(u))| , |e^{z \gt^N(\vl_i(u'))} h(\vl_i(u'))| \}}\\
& =    & \frac{|e^{z \gt^N(\vl_i(u))} h(\vl_i(u)) 
- e^{z \gt^N(\vl_i(u'))} h(\vl_i(u'))|}{ e^{c \gt^N(\vl_i(u'))} |h(\vl_i(u'))| }\\
& \leq & \frac{|e^{z \gt^N(\vl_i(u))} - e^{z \gt^N(\vl_i(u'))} |}{e^{c \gt^N(\vl_i(u'))}}
+ \frac{e^{c \gt^N(\vl_i(u))} | h(\vl_i(u)) - h(\vl_i(u'))|}{ e^{c \gt^N(\vl_i(u'))} |h(\vl_i(u'))| }  \\
& \leq & \frac{|e^{z \gt^N(\vl_i(u))} - e^{z \gt^N(\vl_i(u'))} |}{e^{c \gt^N(\vl_i(u'))}}
+ \frac{ e^{c (\gt^N(\vl_i(u')) - \gt^N(\vl_i(u')))}\,E |b| H(\vl_i(u'))}{ |h(\vl_i(u'))| } D (\vl_i(u),\vl_i(u')) \\
& \leq & \frac{|e^{c \gt^N(\vl_i(u))} - e^{c \gt^N(\vl_i(u'))} |}{e^{c \gt^N(\vl_i(u'))}}
+|e^{\i w \gt^N(\vl_i(u))} - e^{\i w \gt^N(\vl_i(u'))} |
+  4 E|b| {e^{2 a_0 N T}}\, \diam(\cc')\\
& \leq & (e^{C_1t} C_1 t + |w| C_1 t) \, D (\vl_i(u),\vl_i(u')) + 4 E|b| e^{2N a_0 T}\,\frac{\gamma^{n_1} \ep_1}{c_0\gamma^{N}} \\
& \leq & \frac{(B+A_0)\gamma^{n_1} \ep_1}{c_0\gamma^N} + \frac{4 E \gamma^{n_1} \ep_1}{c_0 (e^{-2a_0 T} \gamma)^{N}}  < \frac{\pi}{12}
\end{eqnarray*}
assuming $a_0 > 0$ is is chosen sufficiently small and $N$ sufficiently large.
So, the angle between the complex numbers 
$$e^{z \gt^N(\vl_i(u)} h(\vl_i(u)) \quad  \mbox{\rm and } e^{z \gt^N(\vl_i(u')} h(\vl_i(u'))$$ 
(regarded as vectors in $\R^2$)  is  $< \pi/6$.  In particular,  for each $i = 1,2$ we can choose a real continuous
function $\theta_i(u)$, $u \in  \hZ$, with values in $[0,\pi/6]$ and a constant $\lambda_i$ such that
$$\di e^{z \gt^N(\vl_i(u))} h(\vl_i(u)) = e^{\i(\lambda_i + \theta_i(u))} e^{c \gt^N(\vl_i(u))} |h(\vl_i(u))|$$ 
for all $u\in \hZ$. Fix an arbitrary $u_0\in \hZ$ and set $\lambda = \gamma_\ell(u_0)$. 
Replacing e.g $\lambda_2$ by $\lambda_2 +  2m\pi$ for some integer $m$, we may assume that 
$|\lambda_2 - \lambda_1 + \lambda| \leq \pi$.
Using the above, $\theta \leq 2 \sin \theta$ for $\theta \in [0,\pi/6]$, and some elementary geometry  yields
$|\theta_i(u) - \theta_i(u')|\leq 2 \sin |\theta_i(u) - \theta_i(u')| < \frac{c_2\ep_1}{8}.$

The difference between the arguments of the complex numbers
$$e^{\i \,b\,\tau^N(\vl_1(u))}e^{z \gt^N(\vl_1(u)}  h(\vl_1(u)) \quad \mbox{\rm and} 
\quad e^{\i \,b\, \tau^N(\vl_2(u))}e^{z \gt^N(\vl_2(u)}  h(\vl_2(u))$$
is given by the function
$$\Gl(u) = [b\,\tau^N(\vl_2(u)) + \theta_2(u) + \lambda_2] -  [b\, \tau^N(\vl_1(u)) + \theta_1(u) + \lambda_1]
=  (\lambda_2-\lambda_1) + \gamma_\ell(u) + (\theta_2(u) - \theta_1(u))\;.$$
Given $u'\in \hZ_{j'}$ and $u''\in \hZ_{j''}$, since $\hdd_{j'}$ and $\hdd_{j''}$ are 
contained in $\cc_m$ and are $\eta_\ell$-separable in $\cc_m$, it follows from  Lemma 9 and the above that
\begin{eqnarray*}
|\Gl(u')- \Gl(u'')| \geq  |\gl(u') - \gl(u'')| - |\theta_1(u')-\theta_1(u'')| 
- |\theta_2(u')-\theta_2(u'')| \geq   \frac{c_2\ep_1}{2}\;.
\end{eqnarray*}
Thus,  $|\Gl(u')- \Gl(u'')|\geq \frac{c_2}{2} \epsilon_{1}$ for all $u'\in \hZ_{j'}$ and  $u''\in \hZ_{j''}$. Hence either 
$|\Gl(u')| \geq \frac{c_2}{4}\epsilon_{1}$ for all $u'\in \hZ_{j'}$ or  $|\Gl(u'')| \geq \frac{c_2}{4}\epsilon_{1}$ 
for all $u''\in \hZ_{j''}$.

Assume for example that $|\Gl(u)| \geq \frac{c_2}{4}\epsilon_{1}$ for all 
$u\in \hZ_{j'}$. Since $\hZ \subset \sigma^{n_1}(\cc_m)$, as in \cite{St2} we have
for any $u \in \hZ$ we get $|\Gamma_\ell(u)| < \frac{3\pi}{2}$.
Thus, $\frac{c_2}{4}\epsilon_{1} \leq |\Gl(u)| <  \frac{3\pi}{2}$ for all $u \in \hZ_{j'}$. Now
as in \cite{D} (see also \cite{St2}) one shows that $\chi_\ell^{(1)} (u)  \leq 1$ and $\chi_\ell^{(2)} (u)  \leq 1$ 
for all $u \in \hZ_{j'}$. 
 \hfill\qed

\bs

Parts (a) and (b) of the following lemma can be proved in the same way as the corresponding parts of of Lemma 5.3 in \cite{St2},
while part (c) follows from Lemma 5(b).

\begin{lem}
There exist a positive integer $N$ and constants $\hrho = \hrho(N) \in (0,1)$, $a_0 = a_0(N) > 0$,
$b_0 = b_0(N) > 0$ and $E \geq 1$ such that for every $a, b,c, t, w\in \R$ with $|a|, |c|\leq a_0$, $|b| \geq b_0$
such that (4.1) hold, there exists a finite family $\{ \nn_J\}_{J\in \J}$ of  operators 
$$\nn_J  = \nn_J(a,b,t,c) : \clip_D (\hU) \longrightarrow \clip_D (\hU) ,$$
where $\J = \J(a,b,t,c)$, with the following properties:

(a) {\it The operators $\nn_J$ preserve the cone} $K_{E|b|} (\hU)$ ;

(b) {\it For all $H\in K_{E|b|}(\hU)$ and $J \in \J$ we have}
$\di \int_{ \hU} (\nn_J H )^2 \; d\nu_0 \leq \hrho \; \int_{\hU} H^2 \; d\nu_0$.

(c) {\it If $h, H\in \clip_D (\hU)$ are such that $H\in K_{E|b|}(\hU)$, $|h(u)|\leq H(u)$ for all   $u \in \hU$ and \\
$| h(u) - h(u') | \leq E|b| H(u')\, D (u,u')$ whenever $u,u'\in \hU_i $ for some $i = 1, \ldots,k$,  
then there exists  $J\in \J$  such that $|\labw^N h(u)|\leq (\nn_J H)(u)$ for all $u \in \hU$   and for $z= c + \i w$ we have
$$|(\labtz^N h)(u) - (\labtz^N h)(u')|\leq E|b| (\nn_J H)(u')\, D (u,u')$$
whenever $u,u'\in \hU_i $ for some $i = 1, \ldots,k$.}
\end{lem}

\ms

\noindent
{\it Proof of Theorem} 5(a). 
Using an argument from \cite{D} one derives from Lemma 12 that there exist a positive integer $N$ and constants 
$\hrho \in (0,1)$ and  $a_0 > 0$, $b_0 \geq 1$, $A_0 > 0$ such that for any $a, b, c,t,w\in \R$ with $|a|, |c| \leq a_0$, 
$|b|\geq  b_0$ for which (4.1) hold, and for any $h\in \clip(\hU)$ with $\|h \|_{\lip,b} \leq 1$ we have
\be
\int_{U} |\labtz^{N m}h|^2\; d\nu_0 \leq \hrho^m \quad, \quad m \geq 0 .
\ee

Then the estimate claimed in Theorem 5(a) follows as in \cite{D} (see also the proof of  Corollary 3.3(a) in \cite{St1}).
\hfill\qed

\bs

The proof of Theorem 5(b) can be derived using an approximation procedure as in \cite{D} -- see the Appendix below for
some details.


\def\tfa{\tilde{f}_a}
\def\tmac{ \tilde{\mm}_{ac}}

\def\hmu{\hat{\mu}}

\section{Spectral estimates when $w$ is the leading parameter}
\renewcommand{\theequation}{\arabic{section}.\arabic{equation}}
\setcounter{equation}{0}

Here we try to repeat the arguments from the previous section however changing the roles of
the parameters $b$ and $w$. We continue to use the assumptions made at the beginning of Sect. 4,
however now we suppose that $f \in \clip(\hU)$. We will consider the case
\be
|b| \leq B\,  |w|
\ee
for an arbitrarily large (but fixed) constant $B > 0$. 

Assume that $G : \mt \longrightarrow \R$ is a Lipschitz functions which is constant on stable leaves of 
$B_i = \{ \phi_t(x) : x\in R_i, 0 \leq t \leq \tau(x)\}$
for each rectangle $R_i$ of the Markov family and $A = \min_{x\in \mt} G(x) > 0$. Set
$$L = \Lip(G) \quad ,   \quad D = \diam(\mt)\,,$$
where without loss of generality we may assume that $D \geq 1$. We will also assume that
\be
L \leq \hmu\, A \quad , \mbox{\rm where } \quad \hmu = \frac{c_0\, \hd}{128\, C_0\, C_1\, D} .
\ee
The function
$$g(x) = \int_0^{\tau(x)} G(\phi_t(x))\, dt  \quad, \quad x\in R,$$
is constant on stable leaves of $R$, so it can be regarded as a function on $U$. Clearly $g \in \clip(\hU)$.

\ms

\noindent
{\bf Remark.} Notice that if we replace $G$ by $G+d$ for some constant $d > 0$, then 
$$g'(x) = \int_0^{\tau(x)} (G(\phi_t(x))+d)\, dt = g(x) + d\, \tau(x) ,$$ 
so 
$$\lc_{\fa - \i\,b\tau + \i\, w g} = \lc_{\fa - \i\,b\tau + \i\, w (g'- d\tau)} = \lc_{\fa - \i\,(b + d w)\tau - \i\, w g'} .$$
Choose and fix $d > 0$ so that $\frac{\Lipf(G)}{G_0 + d} \leq \hmu$.
Then for $G' = G+d$ and $g' = g + d \tau$ we have $\frac{\Lipf(G')}{\min G'} \leq \hmu$,
and the operator $\lc_{\fa - \i\,b\tau + \i\, w g} = \lc_{\fa - \i\,b'\tau + \i\, w g'}$,
where $b' = b + dw$. Thus, without loss of generality we may assume that $\frac{\Lip(G)}{\min G} \leq \hmu$,
which is equivalent to (5.2). As in \cite{PeS}, this will imply a non-integrability property for $g$ (see Lemma 10 below).
In other words, dealing with an initial function $G$ one has to first change it to arrange (5.2), and then
with the new parameters $b$ and $w$ that appear in front of $\i \tau$ and $\i g$ consider the cases
$|w| \leq B |b|$ (as in Theorem 5(c)) and $|b|\leq B |w|$, which is considered in this section.

\ms

As in Sect. 4, we will use the set-up and some arguments from \cite{St2}. 
Let $\rr = \{R_i\}_{i=1}^k$ be a Markov family for $\phi_t$ over $\mt$ as in Sect. 1.

Here we prove the following analogue of Theorem 5(c).

\begin{thm} 
Let $\phi_t : M \longrightarrow M$ be a $C^2$ flow satisfying the Standing Assumptions over the basic set $\mt$.
Assume in addition that {\rm (5.2)} holds. Then 
for any real-valued functions $f,g \in \clip(\hU)$,
any constants $\epsilon > 0$ and $B > 0$  there exist constants $0 < \rho < 1$, 
$a_0 > 0$,  $w_0 \geq 1$ and  $C = C(B, \epsilon)> 0$ such that if $a, c\in \R$  satisfy $|a|, |c| \leq a_0$, 
then
\begin{equation} \label{eq:est}
\|L_{f -(P_f+a+ \i b)\tau + (c+ \i w) g}^m h \|_{\lip,b} 
\leq C \, \rho^m \, |b|^{\ep}\, \| h\|_{\lip,b}
\end{equation}
for all integers $m \geq 1$ and all $b,w\in \R$ with  $|w| \geq w_0$ and $|b| \leq B \, |w|$.
\end{thm}

\medskip

Recall the definitions of $\lambda_0 > 0$, $\hnu_0$, $h_{0}$, $\f0$ from Sect. 3; now we have
$h_0, \f0 \in \clip(\hU)$. Fix a small $a_0 > 0$.
Given a real number $a$ with $|a| \leq a_0$, denote by $\lambda_{a}$ the 
{\it largest eigenvalue} of $L_{f -(P+a)\tau}$ on $\clip(U)$ and by $h_{a}$ the corresponding (positive) 
eigenfunction such that $\int_U h_{a}\, d\nu_{a} = 1$, where $\nu_{a}$ is the unique probability measure on 
$U$ with $L_{f -(P+a)\tau}^*\nu_{a} = \nu_{a}$. 
Given real numbers $a, b,c,w$ with $|a|, |c| \leq a_0$ consider the function
$$\tfa = f - (P+a) \tau + \ln h_{a} - \ln (h_{a}\circ \sigma) - \ln \lambda_{a}$$
and the operators 
$$ \labz = L_{\tfa - \i\, b\, \tau + z g} : C(U) \longrightarrow C(U)\:\:\: , \:\:\:
\tmac = L_{\tfa + c g} : C(U) \longrightarrow C(U) ,$$
where $z = c + \i w$. Notice that $L_{\tfa} 1 = 1$.

Taking the constant $C_0 > 0$ sufficiently large, we may assume that
\be
\Lip(\tfa - \f0) \leq C_0 |a| \quad , \quad , \|\tfa - \f0\|_0 \leq C_0\, |a|  \quad, \quad |a| \leq a_0 .
\ee
Thus, ssuming $a_0 > 0$ is chosen sufficiently small, there exists a constant $T > 0$ 
(depending on $f$ and $a_0$) such that 
\be
\|\tfa\|_\infty \leq T  \quad , \quad\Lip(h_{a}) \leq T \quad , \quad \Lip(\tfa) \leq T
\ee
for $|a| \leq a_0$. As before, we will assume that $T \geq \max \{ \, \|\tau \|_0 \, , \, \Lip(\tau_{|\hU}) \, \}$,
and also that $\Lip(g) \leq T$ and $\|g\|_0 \leq T$.

\ms

Essentially in what follows we will repeat (a simplified version of) the proof of Theorem 5, so we will use the set-up in Sect. 4
-- see the text after Lemma 6, up to and including the definition of $\ep_1$.

Let $a,b,c,w\in \R$ be so that $|a|, |c| \leq a_0$, $|w| \geq w_0$, where $w_0$ is a sufficiently large constant defined as $b_0$ 
in Sect. 4, and $|b| \leq B |w|$. Set $z = c + \i w$.

Let $\cc_m$  ($1\leq m \leq p$) be the family of {\it maximal  closed cylinders} 
contained in $\overline{U_0}$ with  $\diam(\cc_m )\leq \frac{\ep_1}{|w|}$ such that 
$U_0 \subset \cup_{j=m}^p \cc_m$ and $\overline{U_0} = \cup_{m=1}^p \cc_m$.  As before we have
$$ \rho\, \frac{\ep_1}{|w|} \leq \diam(\cc_m ) \leq  \frac{\ep_1}{|w|} \quad, \quad 1\leq m\leq p .$$
Fix an integer $q_0 \geq 1$ as in Sect. 4, and let
$\dd_1, \ldots, \dd_q$  be the list of all closed cylinders contained in $\overline{U_0}$ 
that are {\it subcylinders of co-length} $p_0\, q_0$ of some $\cc_m$ ($1\leq m\leq p$). 
Then  $\overline{U_0} = \cc_1 \cup \ldots \cup \cc_p =  \dd_1 \cup \ldots \cup \dd_q$ and
$$\rho^{p_0\, q_0+1}\cdot \frac{\ep_1}{|w|} \leq \diam(\dd_j ) \leq  \rho^{q_0}\cdot 
\frac{\ep_1}{|w|} \quad, \quad 1\leq j\leq q .$$
Next, define the cylinders $Z_j = \overline{\sigma^{n_1}(\hdd_j)}$ and $\xijl = \overline{\vl_i(\hZ_j)}$ as in Sect. 4, and consider 
the characteristic functions  $\omega_{i,j}^{(\ell)} = \chi_{\hxijl} : \hU \longrightarrow [0,1]$.
Let $J$  be a {\it subset of the set} 
$\Xi = \Xi(a,w) = \{\; (i,j, \ell) \; :  \; 1\leq i \leq 2\; ,\;  1\leq j\leq q\; , \;  1\leq \ell \leq \ell_0\;\}.$
Define $\mu_0 > 0$ as in Sect. 4 and $\omega = \omega_{J} : \hU \longrightarrow [0,1]$ by
$\di \omega  = 1- \mu_0 \,\sum_{(i, j, \ell) \in J} \eijl$.
Finally define $\nn = \nn_J(a,b,c) : \clip_D (\hU) \longrightarrow \clip_D (\hU)$ by 
$\left(\nn h\right) = \tmac^N (\omega_J \cdot h)$.

Then we have the following analogue of Lemma 5.

\begin{lem} Assuming $a_0 > 0$ is chosen sufficiently small, there exists a constant $A_0 > 0$ 
such that for all $a,c\in \R$ with $|a|, |c|\leq a_0$  the following hold:

\ms

(a)  {\it If $H \in K_E(\hU)$ for some $E > 0$, then 
$$\frac{|(\tmac^m H)(u) - (\tmac^m H)(u')|}{(\tmac^m H)(u')} \leq
A_0 \, \left[ \frac{E}{\gamma_0^m} + 1 \right]\, D (u,u')$$
for all $m \geq 1$ and all $u,u'\in U_i$, $i = 1, \ldots, k$.}

\ms

(b) {\it If the functions $h$ and  $H$ on $\hU$  and $E > 0$  are such that $H > 0$ on $\hU$ and 
$|h(v) - h(v')| \leq E\, H(v')\, D (v,v')$ for any $v,v'\in \hU_i$, $i = 1, \ldots,k$, 
then for any integer $m \geq 1$ and any $b, w\in \R$ with  $|b|, |w|\geq 1$, for $z = c + \i w$  we have 
$$| (\labw^N  h)(u) - (\labw^N h)(u')| \leq E |w| (\nn H)(u')\, D (u,u') .$$
whenever $u,u'\in \hU_i$ for some $i = 1, \ldots,k$. }
\end{lem}

The proof is a simplified version of that of Lemma 5 and we omit it.


\bs

Next, changing appropriately the definition of a dense subset $J$of $\Xi$, Lemma 8 holds again
replacing $K_{E|b|}(\hU)$ by $K_{E|w|}(\hU)$.

Assume that $h, H\in \clip_D (\hU)$ are such that 
\begin{equation}
H\in K_{E|w|}(\hU)\quad , \quad |h(u)|\leq H(u) \:\:\:\:, \:\:\; u\in \hU  ,
\ee
and
\be
|h(u) - h(u')|\leq E|w| H(u')\, D (u,u')\:\:\:\:\:
\mbox{\rm whenever} \:\: u,u'\in \hU_i\;, \; i = 1, \ldots,k  .
\end{equation}  

Define the functions $\geil : \hU  \longrightarrow \C$ by
$$\di \geol(u) = \frac{\di \left| e^{(\tfa^{N} - \i b\tau^{N} + z g^N)(\vl_1(u))} h(\vl_1(u)) +
 e^{(\tfa^{N} - \i b\tau^{N} + z g^N)(\vl_2(u))} h(\vl_2(u))\right|}{\di (1-\mu)
e^{\tfa^{N}(\vl_1(u)) + c g^N( \vl_1(u))}H(\vl_1(u)) + e^{\tfa^{N}(\vl_2(u))+ c g^N(\vl_2(u))}H(\vl_2(u))} ,$$ 
$$\di \getl(u) = \frac{\di \left| e^{(\tfa^{N} - \i b\tau^{N} + z g^N)(\vl_1(u))} h(\vl_1(u)) +
 e^{(\tfa^{N} - \i b\tau^{N} + z g^N)(\vl_2(u))} h(\vl_2(u))\right|}{\di
e^{\tfa^{N}(\vl_1(u)) + c g^N(\vl_1(u))}H(\vl_1(u)) +  (1-\mu) e^{\tfa^{N}(\vl_2(u)) + cg^N(\vl_2(u))}H(\vl_2(u))} ,$$ 
and set
$ \gl(u) = w\, [\tau_N(\vl_2(u)) - \tau_N(\vl_1(u))]$, $u\in \hU$.
The crucial step in this section is to prove the following analogue of Lemma 9:

\begin{lem}
Let $j, j'\in \{ 1, 2,\ldots,q\}$ be
such that $\dd_j$ and $\dd_{j'}$ are contained in $\cc_m$ and  are $\eta_\ell$-separable 
in $\cc_m$ for some $m = 1, \ldots, p$ and $\ell = 1, \ldots, \ell_0$ . Then
$ |\gl (u) - \gl(u')| \geq c_3\epsilon_{1}$ for all $u\in \hZ_j$ and  $u'\in \hZ_{j'}$,
where $\di c_3 =  \frac{A \hd\, \rho}{32}$.
\end{lem}

To prove the above we need the following.

\begin{lem} 
{\rm (Lemma 6 in \cite{PeS})}
Assume that {\rm (5.2)} holds. 
Under the assumptions and notation in Lemma 1, for all  $\ell = 1, \ldots, \ell_0$, 
$s\in r^{-1}(U_0)$, $0 < |h| \leq \hd$ and $\eta \in B_\ell$ so that  
$s+h\, \eta \in r^{-1}(U_0\cap \mt)$ we have 
$$\left[I_{\eta,h} \left(g^{N}(\vl_2(\trr(\cdot ))) - g^{N}(\vl_1(\trr(\cdot)))\right)\right](s)  \geq \frac{A \hd}{4}\,.$$
\end{lem}

\ms

\noindent
{\it Proof of Lemma} 14. This just a repetition of the proof of Lemma 5.9 in \cite{St2}, where instead of using
Lemma 6(b) we use the above Lemma 14. We omit the details.
\hfill\qed

\bs

Next, we need to prove the analogue of Lemma 10.

\begin{lem}
Assume $|w| \geq w_0$ for some sufficiently large $w_0 > 0$ and let $|b| \leq B |w|$. Then  for any
$j = 1, \ldots,q$  there exist $i \in \{ 1,2\}$, $j' \in \{ 1,\ldots,q\}$ and $\ell \in \{ 1, \ldots, \ell_0\}$  
such that  $\dd_j$ and $\dd_{j'}$ are adjacent and $\chi^{(i)}_{\ell} (u) \leq 1$ for all  $u\in \hZ_{j'}$. 
\end{lem}

\ms

\noindent
{\it Sketch of proof of Lemma } 16. We will use Lemma 11 which holds again with (4.3)-(4.4) replaced by (5.6)-(5.7).

Given $j = 1, \ldots, q$, let $m = 1, \ldots, p$ be such that $\dd_j \subset \cc_m$.  
As in \cite{St2} we find $j',j'' = 1, \ldots,q$ such that $\dd_{j'}, \dd_{j''} \subset \cc_m$ and 
$\dd_{j'}$ and $\dd_{j''}$ are $\eta_\ell$-separable in $\cc_m$. 

Fix $\ell$, $j'$ and $j''$ with the above properties, and set $\hZ = \hZ_j \cup \hZ_{j'}\cup \hZ_{j''}\; .$ 
If there exist $t \in \{j, j', j''\}$ and $i = 1,2$ such that the first alternative in Lemma 11(b) holds for 
$\hZ_{t}$, $\ell$  and $i$, 
then $\mu \leq 1/4$ implies $\chi_\ell^{(i)}(u) \leq 1$ for any $u\in \hZ_{t}$.

Assume that for every $t\in \{ j, j', j''\}$ and every $i = 1,2$ the second alternative 
in Lemma 11(b) holds for $\hZ_{t}$, $\ell$ and $i$, i.e. 
$|h(\vl_i(u))|\geq \frac{1}{4}\, H(\vl_i(u))$, $u \in \hZ$.

Again we have $\psi(\hZ) = \hdd_j \cup \hdd_{j'} \cup \hdd_{j''} \subset \cc_m$, and 
$\cc' = \vl_i (\sigma^{n_1}(\cc_m))$ is a cylinder with $\diam(\cc') \leq \frac{\ep_1}{c_0\, \gamma^{N-n_1}\, |w|}$. 
Thus, 
assuming e.g.   $|h(\vl_i(u))| \geq |h(\vl_i(u'))|$, we get
\begin{eqnarray*}
&      & \frac{|e^{\i b \tau_N(\vl_i(u)} h(\vl_i(u)) 
- e^{\i b \tau_N(\vl_i(u')} h(\vl_i(u'))|}{\min\{ |h(\vl_i(u))| , |h(\vl_i(u'))| \}}\\
& \leq & |e^{\i b \tau_N(\vl_i(u)} - e^{\i b \tau_N(\vl_i(u')} | 
+ \frac{E|w|\, H(\vl_i(u'))}{|h(\vl_i(u'))| } D (\vl_i(u),\vl_i(u')) \\
& \leq & |b| C_1 \, D (\vl_i(u),\vl_i(u')) + 4 E|w|\, D (\vl_i(u),\vl_i(u')) \\
& \leq & (B |w| \,C_1     + 4 E|w|) \, \diam (\cc')
\leq \frac{(B C_1   + 4 E) \ep_1}{\gamma_1^{N-n_1}}  < \frac{\pi}{12}
\end{eqnarray*}
assuming $N$ is chosen sufficiently large.
So, the angle between the complex numbers 
$$e^{\i b \tau_N(\vl_i(u)} h(\vl_i(u)) \quad  \mbox{\rm and } e^{\i b \tau_N(\vl_i(u')} h(\vl_i(u'))$$ 
(regarded as vectors in $\R^2$)  is  $< \pi/6$.  
In particular,  for each $i = 1,2$ we can choose a real continuous
function $\theta_i(u)$, $u \in  \hZ$, with values in $[0,\pi/6]$ 
and a constant $\lambda_i$ such that
$\di h(\vl_i(u)) = e^{\i(\lambda_i + \theta_i(u))}|h(\vl_i(u))|$ for all $u\in \hZ$.
Fix an arbitrary $u_0\in \hZ$ and set $\lambda = \gamma_\ell(u_0)$. 
Replacing e.g $\lambda_2$ by $\lambda_2 +  2m\pi$ for some integer $m$, 
we may assume that $|\lambda_2 - \lambda_1 + \lambda| \leq \pi$.
Using the above, $\theta \leq 2 \sin \theta$ for $\theta \in [0,\pi/6]$, and some elementary geometry  yields
$|\theta_i(u) - \theta_i(u')|\leq 2 \sin |\theta_i(u) - \theta_i(u')| < \frac{c_2\ep_1}{8}.$

The difference between the arguments of the complex numbers
$$e^{\i \,b\,\tau_N(\vl_1(u))}e^{\i w g_N(\vl_1(u)}  h(\vl_1(u)) \quad \mbox{\rm and} 
\quad e^{\i \,b\, \tau_N(\vl_2(u))}e^{\i w g_N(\vl_2(u)}  h(\vl_2(u))$$
is given by the function
$$\Gl(u) = [w\, g_N(\vl_2(u)) + \theta_2(u) + \lambda_2] -  [w\, g_N(\vl_1(u)) + \theta_1(u) + \lambda_1]
=  (\lambda_2-\lambda_1) + \gamma_\ell(u) + (\theta_2(u) - \theta_1(u))\;.$$
Given $u'\in \hZ_{j'}$ and $u''\in \hZ_{j''}$, since $\hdd_{j'}$ and $\hdd_{j''}$ are 
contained in $\cc_m$ and are $\eta_\ell$-separable in $\cc_m$, it follows from  Lemma 9 and the above that
\begin{eqnarray*}
|\Gl(u')- \Gl(u'')| \geq  |\gl(u') - \gl(u'')| - |\theta_1(u')-\theta_1(u'')| 
- |\theta_2(u')-\theta_2(u'')| \geq   \frac{c_3\ep_1}{2}\;.
\end{eqnarray*}
Thus,  $|\Gl(u')- \Gl(u'')|\geq \frac{c_3}{2} \epsilon_{1}$ for all $u'\in \hZ_{j'}$ and 
$u''\in \hZ_{j''}$. Hence either 
$|\Gl(u')| \geq \frac{c_3}{4}\epsilon_{1}$ for all $u'\in \hZ_{j'}$ or 
$|\Gl(u'')| \geq \frac{c_3}{4}\epsilon_{1}$ for all $u''\in \hZ_{j''}$.

Assume for example that $|\Gl(u)| \geq \frac{c_2}{4}\epsilon_{1}$ for all 
$u\in \hZ_{j'}$. Since $\hZ \subset \sigma^{n_1}(\cc_m)$, as in \cite{St2} we have
for any $u \in \hZ$ we get $|\Gamma_\ell(u)| < \frac{3\pi}{2}$.
Thus, $\frac{c_2}{4}\epsilon_{1} \leq |\Gl(u)| <  \frac{3\pi}{2}$ for all $u \in \hZ_{j'}$. Now
as in \cite{D} (see also \cite{St2}) one shows that
$\chi_\ell^{(1)} (u)  \leq 1$ and $\chi_\ell^{(2)} (u)  \leq 1$ for all $u \in \hZ_{j'}$.  
\hfill\qed

\bs

\noindent
{\it Proof of Theorem} 6. This is now the same as the proof of Theorem 5(a).
\hfill\qed

\section{Analytic continuation of the function $\zeta(s, z)$}
\renewcommand{\theequation}{\arabic{section}.\arabic{equation}}
\setcounter{equation}{0}

Consider the function $\zeta(s, z)$ introduced in Section 1. Recall that $s = a + \i b, z = c + \i w$ with real 
$a, b, c, w \in \R$.  First, we assume that $f$ and $g$ are functions in $C^{\alpha}(\Lambda)$ with some $0 < \alpha < 1$. 
Passing to the symbolic model defined by the Markov family ${\mathcal R}$ we obtain functions\footnote{In fact, one has to define
first $f$ and $g$ as functions in $C^{\alpha}(\hat{R})$ and then extend them as $\alpha$-H\"older functions on $R$. 
In the same way one should proceed with H\"older functions on $U$.} in $C^{\alpha}(R)$ which we denote again by $f$ and $g$. 
We assume that $Pr(f - P_f \tau) = 0$ and we set $s = P_f + a + \i b$. The functions $f$, $g$  depend on $x \in R$. 
A second reduction is to replace $f$ and $g$ by functions $\hat{f}, \: \hat{g} \in C^{\alpha/2}(U)$ depending only on 
$x \in U$ so that $f = \hat{f} + h_1 - h_1 \circ \sigma, \: g = \hat{g} + h_2 - h_2 \circ \sigma$ (see Proposition 1.2 in \cite{PP}). 
Since for periodic points with $\sigma^nx = x$ we have $f^n(x) = \hat{f}^n(x),\: g^n(x) = \hat{g}^n (x)$, we obtain the 
representation
$$\zeta(s, z) = \exp\Bigl( \sum_{n=1}^{\infty} \frac{1}{n} \sum_{\sigma^n x = x}  
e^{\hat{f}^n(x) -(P_f + a + \i b) \tau^n(x)+ (c + \i w) \hat{g}^n(x)}\Bigr).$$ 

In this section we will prove under the standing assumptions that there exists $\epsilon > 0$ and $\epsilon_0 > 0$ such that 
the function $\zeta(s, z)$ has a non non zero analytic continuation for
$ - \epsilon \leq a \leq 0$ and  $|z| \leq \epsilon_0$ with a simple pole at $s = s(z),\: s(0) = P_f.$ Here $s(z)$ is determined 
from the equation $Pr(f - s\tau + z g) = 0.$ For simplicity of the notation we denote below $\hat{f}$ and $\hat{g}$ again by 
$f$, $g$.\\

First consider the case $0 < \delta \leq |b| \leq b_0.$ Since our standing assumptions imply that the flow
$\phi_t$ is weak mixing, Theorem 6.4 in \cite{PP} says that for every fixed $b$ lying in the compact interval 
$[\delta, b_0]$ there exists $\epsilon(b) > 0$ so that the function $\zeta(s, z)$ is analytic for 
$| s - P_f + \i b | \leq \epsilon(b), \: |z| \leq \epsilon(b).$ This implies that there exists $\eta_0 = \eta_0(\delta, b_0) > 0$ 
such that $\zeta(s, z)$ is analytic for 
$P_f - \eta_0 \leq \Re s \leq P_f + \eta_0, \delta \leq |\Im s| \leq b_0,\: |z| \leq \eta_0.$ Decreasing $\delta > 0$ and 
$\eta_0$, if it is necessary, we apply once more Theorem 6.4 in \cite{PP}, to conclude that 
$\zeta(s, z) ( 1 - e^{Pr(f - s \tau + z g)})$ is analytic for 
$$s \in \{ s \in \C:\:|\Re s -P_f| \leq \eta, |\Im s | \leq \delta\}$$
and $|z| \leq \eta_0.$ Consequently, the singularities of $\zeta(s, z)$ are given by $(s, z)$ for which we have 
$Pr(f - s\tau + z g) = 0$ and, solving this equation, we get $s = s(z)$ with $s(0) = P_f$. It is clear that we have a 
simple pole at $s(z)$ since $\frac{d}{ds} Pr(f - s \tau + z g) \neq 0$ for $|z|$ small enough.\\

Now we pass to the case when $|\Im s| = |b| \geq b_0 > 0,\: |z| \leq \eta_0.$ Then we fix a $\beta \in (0, \alpha/2)$ and we get 
with $0 < \mu < 1$ the inequality $|\Im b| \geq B_0 |z|^{\mu} $ with  $B_0 = \frac{b_0}{\eta_0^{\mu}}$. Thus we are in position 
to apply the estimates of Theorem 5(b) saying that for every $\epsilon > 0$ there exist 
$0 < \rho < 1$ and $C_{\epsilon} > 0$ so that
\begin{equation} \label{eq:6.1}
\|L_{f - (P_f + a + \i b)\tau + z g}^m\|_{\beta,  b} \leq C_{\epsilon} \rho^m |b|^{\epsilon} ,\:\forall m \in \N
\end{equation}
for $|a| \leq a_0, |b| \geq b_0, \: |z| \leq \eta_0.$ Next we apply Theorem 4 with functions $f, g \in C^{\beta}(U)$. 
For $|\Re s - P_f| \leq \eta_0,\:|\Im s | \geq b_0$  and $|z| \leq \eta_0$ we deduce 
 $$|Z_n(f - (P_f + a + \i b) \tau + z g)| \leq \sum_{i = 1}^k |L^n_{f - (P_f + \i a + b)\tau + z g} (\chi_i)(x_i)|$$
$$+ C( 1 + |b|) \sum_{m = 2}^n \|L_{f -(P_f + a + \i b) \tau + z g}^m\|_{\beta} \gamma_0^{-m \beta}e^{m Pr(f-(P_f + a) \tau 
+ (\Re z) g)}$$ 
$$\leq k \|L_{f - (P_f + a + \i b) \tau + z g}^n\|_{\beta} + C_{\epsilon}(1 
+ |b)|b|^{\epsilon} \sum_{m = 2}^n \rho^{n - m} \gamma_0^{-m \beta} e^{m(\epsilon +  Pr(f - ( P_f + a )\tau + c g))}.$$
Taking $\eta_0$  and $\epsilon$ small, we arrange
$$\gamma_0^{-\beta} e^{\epsilon +  Pr(f - (P_f + a)\tau + c g)} \leq \gamma_2 < 1$$
for $|a| \leq \eta_0, \: |c| \leq \eta_0$, since $Pr(f - P_f \tau) = 0$ and $\gamma_0^{-\nu} < 1.$ Next increasing $0 < \rho < 1$, 
if it is necessary, we get $\frac{\gamma_2}{\rho} < 1.$ Thus the sum above will be bounded by
$$C_{\epsilon}(1 + |b|)|b|^{\epsilon} \rho^n \sum_{m= 2}^{\infty} \Bigl(\frac{\gamma_2}{\rho}\Bigr)^m 
\leq C'_{\epsilon} |b|^{1 + \epsilon} \rho^n $$
for $|a| \leq \eta_0, \: |z| \leq \eta_0.$ The analysis of the term $\|L_{f -(P_f + a + \i b) + z g}^n\|_{\beta}$ follows the 
same argument and it is simpler. Finally, we get
$$|Z_n(f -(P_f + a + \i b) \tau + z g)| \leq B_{\epsilon} |b|^{1 + \epsilon}\rho^n,\: \forall n \in \N$$
and the series 
$$\sum_{n=1}^{\infty}\frac{1}{n} Z_n(f -(P_f + a + \i b)\tau + z g)$$
is absolutely convergent for $|a| \leq \eta_0, |b| \geq b_0, |z| \leq \eta_0$. This implies the analytic continuation of 
$\zeta(s, z)$ for $|\Re s - P_f| \leq \eta_0, |\Im s| \geq b_0, \: |z| \leq \eta_0$,
thus completing the proof of Theorem 1.\\

To obtain a representation of the function $\eta_g(s) = \frac{\partial\log\zeta(s, z)}{\partial z} \big\vert_{z = 0}$ for $s$ sufficiently close to $P_f$, notice that for such values of $s$ we have
$$\eta_g(s) = - \frac{\partial\log(1 - e^{Pr(f - s\tau + zg)})}{\partial z}\big\vert_{z = 0} + A_0(s)$$
$$= \frac{1}{s -P_f}\frac{\int g d m}{\int \tau dm} + A_1(s) = \frac{\int G d \mu_F}{s - P_f} + A_1(s),$$
where $m$ is the equilibrium state of $f - P_f\tau$, $\mu_F$ is the equilibrium state of $F$ and $A_0(s)$ and $A_1(s)$ are analytic in a neighborhood of $P_f$ (see Chapter 6 in \cite{PP}).  More precisely, $\mu_F$ is a $\sigma_t^{\tau}$ invariant probability measure on $R^{\tau}$ such that
$$Pr(F) = h(\sigma^{\tau}_1, \mu_F) + \int F (\pi(x, t))d\mu_F,$$
where $h(\sigma^{\tau}_1,\mu_F)$ is the metric entropy of $\sigma^{\tau}_1$ with respect to $\mu_F$ (see Chapter 6 in \cite{PP}).

Taking $\eta_0$ small enough,  for $|z| \leq \eta_0,\: |\Re s - P_f| \leq \eta_0$ and  $|\Im s| \geq  \eta_0$ from the estimates for $Z_n(f - (P_f + a + \i b) \tau + z g)$ above, we deduce
$$|\log \zeta(s, z)| \leq C_{\epsilon} \max\Bigl(1, |\Im s|^{1 + \epsilon}\Bigr).$$ 
To estimate $\eta_g(s)$, as in \cite{PoS2}, we apply the Cauchy theorem for the derivative
$$\frac{\partial}{\partial z} \log\zeta(s, z)\big\vert_{z = 0} = \frac{1}{2 \pi \i \delta} \int_{|\xi| = \delta} \frac{\log\zeta(s, \xi)}{\xi^2} d\xi= {\mathcal O}(|\Im s|^{1 + \epsilon}),\: |\Im s | \geq 1.$$
with $\delta > 0$ sufficiently small. Thus we obtain a ${\mathcal O}\Bigl(\max\Bigl(1, |\Im s|^{1 + \epsilon}\Bigr)\Bigr)$ bound for the function 
$$A(s) = \eta_g(s) - \frac{1}{s - P_f} \int G d\mu_F$$
which is analytic for $|\Re s - P_f| \leq \eta_0$. Decreasing $\eta_0$ and applying Phragm\'en-Lindel\"of theorem, by a standard argument we obtain a bound ${\mathcal O}\Bigl(\max\Bigl(1, |\Im s|^{\alpha}\Bigr)\Bigr)$ with $0 < \alpha < 1.$ Consequently, we have the following
\begin{prop} Under the assumptions of Theorem $1$ there exist $\eta_0 > 0$ and $ 0 < \alpha < 1$ such that for 
$\Re s > P_f - \eta_0$ we have
\begin{equation} \label{eq:6.2}
\eta_g(s) =  \frac{1}{s - P_f} \int G d\mu_F + A(s)
\end{equation}
with an analytic function $A(s)$ satisfying the estimate
\begin{equation} \label{eq:6.3}
|A(s)| \leq C \max\Bigl(1, |\Im s|^{\alpha}\Bigr). 
\end{equation}
\end{prop}

Next  define  $\ff^{\tau}(\C):=\{ F: R^{\tau} \longrightarrow \C\}$ and $\ff^{\tau}(\R) :=\{F: R^{\tau} \longrightarrow \R\}$ 
the spaces of complex-valued (real-valued) functions which are continuous. If $G \in \ff^{\tau}(\C)$ is Lipschitz continuous 
and if the standing assumptions for $\Lambda$ are satisfied, 
 the function
$$g(x) = \int_0^{\tau(x)} G(\pi(x, t)) dt$$
is Lipschitz continuous on $R$. Moreover, if the representative of $G$ in the suspension space $R^{\tau}$ is constant on stable leaves, the function $g(x)$ depends only on $x \in U.$ 
Now we introduce two definitions of independence.

\begin{deff} Two functions $f_1, f_2: U \to \R$ are called $\sigma-$ independent if whenever there are constants 
$t_1, t_2 \in R$ such that
$t_1 f_1 + t_2 f_2$ is co homologous to a function in $C(U : 2 \pi \Z)$, we have $t_1 = t_2 = 0.$
\end{deff}

For a function $G \in \ff^{\tau}(\R)$ consider the skew product flow $S_t^G$ on ${\mathbb S}^1 \times R^{\tau}$ by
$$S_t^G ( e^{2 \pi \i \alpha}, y) = \Bigl(e^{2 \pi \i (\alpha + G^t(y))}, \sigma_t^{\tau}(y)\Bigr).$$
\begin{deff} [\cite{La}] Let $G \in \ff^{\tau}(\R)$. Then $G$ and $\sigma^{\tau}_t$ are {\it flow independent} if 
the following condition is satisfied. If $t_0, t_1 \in \R$ are constants such that the skew product flow $S_t^H$ with 
$H = t_0 + t_1 G$ is not topologically mixing, then $t_0 = t_1 = 0.$
\end{deff}
Notice that if $G$ and $\sigma^{\tau}_t$ are flow independent, then the flow $\sigma^{\tau}_t$ is topologically weak mixing 
and the function $G$ is not co homologous to a constant function. On the other hand, if $G$ and $\sigma^{\tau}_t$ are flow
 independent, then $g(x) = \int_0^{\tau(x)} G(\pi(x, t)) dt$ and $\tau$ are $\sigma-$ independent.\\ 
Below we assume that $g$ and $\tau$ are $\sigma-$ independent and we suppose that $F, G$ is a Lipschitz functions $\Lambda$ having representative in $R^{\tau}$ which are constant on stable leaves. Thus we obtain functions $f$, $g$ which are in $C^{\Lip}(\hU)$.   
We will now obtain an analytic continuation of 
$\zeta(s, z)$ for $P_f - \eta_0 < \Re s < P_f$ and $z = \i w$. Set $r(s, w) = f - (P_f + a + \i b)\tau + \i w g.$ 
We choose $M > 0$ large enough so that we can apply Theorem 6 for $|w| \geq M.$ We consider two cases.\\

 {\bf Case 1. $ \eta_0 \leq |w| \leq M.$} We consider two sub cases: 1a) $|\Im s| \leq  M_1, \: 1b)\: |\Im s| \geq  M_1$. 
 Here $M_1 > 0$ is chosen large enough so that Theorem 5 (b) holds with $|\Im s| \geq M_1.$\\

Let $|\Im s| \leq M_1.$ Assume first that $\Im r(s_0, w_0)$ is cohomologous to $c + 2\pi Q$ with an integer-valued function 
$Q \in C(U; \Z)$ and a constant $c \in [0, 2 \pi ).$ If $c = 0$, since $g$ and $\tau$ are $\sigma-$ independent, from the fact 
that $b \tau + w g$ is co homologous to a function in $C(U; 2 \pi \Z),$ we deduce $b = w = 0$ which is impossible because 
$b = \Im s \neq 0.$ Thus we have $c \neq 0.$ Consequently, the operator $L_{f -s_0 \tau + \i wg}$ has an eigenvalue $ e^{\i c}$. 
Then there exists a neighborhood $U_1 \subset \C \times \R$ of $(s_0, w_0)$ such that for $(s, w) \in U_1$ we have 
$Pr(r(s, w)) \neq 0$ and for
$(s, w) \in U_2$ we have an analytic extension of $\log\zeta(s, w)$ given by
$$\log\zeta(s, w) =  \frac{K_1(s, w)}{1 - e^{Pr(r(s, w))}} + J_1(s, w)$$
with functions $K_1(s, w), J_1(s, w)$ analytic with respect to $s$ for $(s, w) \in U_1.$ Second, let $\Im r(s_0, w_0)$ 
be not cohomologous to $c + 2 \pi Q$. Then the spectral radius of $L_{f - s_0 \tau + \i w g}$ is strictly less than 1 
and this will be the case for $(s, w)$ is a small neighborhood $U_2 \subset \C \times \R$ of $(s_0, w_0).$ Applying Theorem 4, 
this implies easily that $\log\zeta(s,\i w)$ has an analytic continuation with respect to $s$.
\\

Passing to the case 1b), we observe that $|\Im s| \geq \frac{M_1}{\eta_0}|w|.$ Then, we apply Theorem 5, (c) combined with Theorem 4 
to obtain an analytic continuation of $\log \zeta(s, \i w).$  Moreover, our argument works for $z = c + \i w$ with 
$|c| \leq \eta_0$ and $\eta_0 \leq |w| \leq M$ and we obtain an analytic continuation of $\log\zeta(s, z)$ for 
$P_f - \eta_0 \leq \Re s < P_f, |c| \leq \eta_0,\: \eta_0 \leq |w| \leq M.$\\

{\bf Case 2. $|w| \geq M$}. We consider two sub cases: 2a) $|\Im s| \geq B |w|,$ 2b) $|\Im s| \leq B|w|, \: B = \frac{M_1}{M}.$ 
If we have 2a), we apply Theorem 5 (c).  In the case 2b) we  use the argument of Section 5 replacing $g(x)$ by 
$g'(x) = g(x) + d \tau(x)$, where the constant $d > 0$ is chosen so that for the function $G'= G + d$ we have
$$\frac{{\rm Lip}\: G'}{\min G'} \leq \hat{\mu},$$
where $\hat{\mu}> 0$ is the constant introduced in Section 5.
Next we write 
$$L_{f - (P_f + a + \i b)\tau + \i w g} = L_{f -(P_f + a + \i (b + d w)\tau + \i w g'}.$$
For the Ruelle operator involving $g'$ we can apply Theorem 6 since $|b + d w|  \leq (B + d)|w|$, $|w| \geq M$ and $g$ is a 
Lipschitz function.
An application of Theorem 4 implies the analytic continuation of $\log \zeta(s, \i w)$ for $P_f - \eta_0 \leq \Re s \leq P_f$ 
and $|w| \geq M.$ From the above analysis we deduce the following

\begin{thm} Assume the standing assumptions fulfilled for the basic set $\Lambda$. Let $F,G: \Lambda \longrightarrow \R$ be 
Lipschitz functions having representatives in $R^{\tau}$ which are constant on stable leaves.
Assume that $g$ and $\tau$ are $\sigma$-independent. Then there exists $\eta_0 > 0$ such that $\zeta(s, \i w)$ admits a non 
zero analytic continuation with respect to $s$ for $P_f - \eta_0 \leq \Re s,\: w \in \R$ and $|w| \geq \eta_0.$
\end{thm}

\section{Applications}
\renewcommand{\theequation}{\arabic{section}.\arabic{equation}}
\setcounter{equation}{0}

\subsection{ Hannay-Ozorio de Almeida sum formula.}

The proof of (\ref{eq:1.5}) in \cite{PoS3} is based on the analytic continuation of the Dirichlet series
$$\eta(s) = \sum_{\gamma} \sum_{m=1}^{\infty} \lambda_G(\gamma) e^{m(-\lambda^u(\gamma) - (s-1) \lambda(\gamma))},\: s \in \C$$
for $1 - \eta_0 \leq \Re s < 1.$ For this purpose the authors examine the analytic continuation of the symbolic function 
$\eta_g(s)$ with $g(x) = \int_0^{\tau(x)} G(\pi(x, t))dt$ defined in Section 1 and they use the fact that the difference 
$\eta(s) - \eta_g(s)$ is analytic in a region $\Re s > 1 - \epsilon',\: \epsilon'> 0.$ Next for the geodesic flow on surfaces 
with negative curvature they establish Proposition 3 with $P_f = 1$. Since $M$ is an attractor, the equilibrium state of the 
function $-E(x)$ is just the SRB measure $\mu$ of $\phi_t$ (see \cite{BR}) and the residuum in (\ref{eq:6.2}) becomes 
$\int G d \mu.$\\
 
For the proof of Proposition 3 in \cite{PoS3} the authors exploit the link between the analytic continuation of $\zeta(s, z)$ 
and the spectral estimates of the Ruelle operator obtained by Dolgopyat \cite{D}. However, in \cite{PoS3} Ruelle's lemma 
in \cite{PoS1} was used whose proof is rather sketchy and contains some steps which are not done in detail (see \cite{W} for more 
information and comments concerning these steps and the gaps in their proofs). On the other hand, the estimates of Dolgopyat 
\cite{D} are established only for Ruelle operators with one complex parameter, and to take into account the second 
parameter $z$ some complementary analysis is necessary.\\

We would like to mention that \cite{W} contains a correct and complete proof of Ruelle's lemma in the case of one complex parameter 
and H\"older function $\tau(x)$. A  version of this lemma with two complex parameters is given in Section 2 above. 
Next, in Theorem 5 the spectral estimates for the  Ruelle operator with two complex parameters are established
for Axiom A flows on a basic set $\Lambda$ of arbitrary dimension under the standing assumptions. If $\Lambda$ is an attractor, 
according to  \cite{BR}, the equilibrium state of $-E(x)$ coincides with the SRB measure $\mu$ of $\phi_t$. Thus we can 
apply Proposition 3 to obtain a representation of $\eta_g(s)$ with residue $\int G d\mu.$ Using (\ref{eq:6.2}) and 
repeating the argument of Section 4 in \cite{PoS3}, we obtain Theorem 2.

\subsection{ Asymptotic of the counting function for period orbits.} As we mentioned in Sect. 1, the analysis 
of $\pi_F(T)$ is based on the analytic continuation of the function  $\zeta(s, 0)$ defined in Section 1. From the arguments
in Section 6 with $z= 0$ and the proof of Proposition 3 we get the following
\begin{prop} Under the standing assumptions in Sect. 3 there exists $\eta_0 > 0$ such that $\frac{\zeta'_F(s)}{\zeta_F(s)}$ admits an analytic
 continuation for $Pr(F) - \eta_0 \leq \Re s$ with a simple pole at $s = Pr(F)$ with residue 1. Moreover, there exists 
 $0 < \alpha < 1$ such that for $|\Im s| \geq 1$ we have
\begin{equation} \label{eq:7.2}
\big|\frac{\zeta'_F(s)}{\zeta_F(s)}\big| \leq C |\Im s|^{\alpha}.
\end{equation}
\end{prop}
To obtain an asymptotic of $\pi_F(T)$, we examine the functions 
$$\Psi(T) = \sum_{e^{n Pr(F)\lambda(\gamma)} \leq T}
\lambda(\gamma) e^{Pr(F) \lambda(\gamma)}, \: \Psi_1(T) = \int_0^T \Psi(y) dy.$$
By a standard argument (see \cite{PoS1} and \cite{Po1}) we obtain the representation
$$\psi_1(T) = \frac{T^2}{2} + \int_{\Re s = (1-\eta_0)Pr(F)}\Bigl(-\frac{\zeta_F'(s)}{\zeta_F(s)}\Bigr) \frac{T^s}{s(s+1)} ds
 = \frac{T^2}{2} + {\mathcal O}(T^{1 + \alpha}),$$ 
where in the second equality the estimate (\ref{eq:7.2}) is used.
This implies an asymptotic for $\Psi(T)$ and repeating the argument in \cite{PoS1}, \cite{Po1}, one obtains Theorem 3.

\bs

\section{Appendix: Proofs of some lemmas}
\renewcommand{\theequation}{\arabic{section}.\arabic{equation}}
\setcounter{equation}{0}

\noindent
{\it Proof of Lemma} 4. Denote by $\ff_\theta(\hU)$ the space of all functions $h : \hU \longrightarrow \R$ that are
Lipschitz with respect to $d_\theta$.  Let $g \in \clip(\hU)$, and let $\theta = \theta_\alpha \in (0,1)$ be as in Sect. 3.
Then $g\in \ff_\theta(\hU)$. Let $\lambda > 0$ be the maximal positive eigenvalue of $L_g$ on $\ff_\theta(\hU)$ and let 
$h > 0$ be a corresponding normalized eigenfunction. By the Ruelle-Perron-Frobenius theorem, we have that
$\frac{1}{\lambda^m} L_g^m 1$ converges uniformly to $h$. We will show that there exists a constant
$C > 0$ such that $\frac{1}{\lambda^m} \Lip(L_g^m 1) \leq C$ for all $m$; this would then imply immediately that
$h\in \clip(\hU)$ and $\Lip(h) \leq C$.

Take an arbitrary constant $K > 0$ such that $1/K \leq h(x) \leq K$ for all $x\in \hU$.
Given $u,u' \in \hU_i$ for some $ i = 1, \ldots,k$ and an integer $m \geq 1$ for any $v\in \hU$ with
$\sigma^m(v) = u$, denote by $v' = v'(v)$ the unique $v'\in \hU$ in the cylinder of length $m$ containing $v$ such
that $\sigma^m(v') = u'$. By (1.1) we have
$$|g_m(v) - g_m(v')| \leq \sum_{j=0}^{m-1} |g(\sigma^j(v)) - g(\sigma^j(v'))|
\leq \Lip(g)\, \sum_{j=0}^{m-1} \frac{d(u,u')}{c_0 \gamma^m} \leq C' \, \Lip(g)\, d(u,u') $$
for some constant $C' > 0$.
Thus,
\begin{eqnarray*}
|(L_g^m1)(u) - (L_g^m 1)(u')|
& \leq & \sum_{\sigma^m(v) = u} \left| e^{g_m(v)} - e^{g_m(v')}\right|
= \sum_{\sigma^m(v) = u} e^{g_m(v)}\, \left| e^{g_m(v)- g_m(v')}-1 \right|\\
& \leq & e^{ C' \, \Lipf(g)}\,  \sum_{\sigma^m(v) = u} e^{g_m(v)}\, \left| g_m(v)- g_m(v') \right|\\
& \leq & e^{ C' \, \Lipf(g)}\,C' \Lip(g) \, d(u,u')\,  \sum_{\sigma^m(v) = u} e^{g_m(v)} \\
& \leq & e^{ C' \, \Lipf(g)}\,C' \Lip(g) \, d(u,u')\,  \sum_{\sigma^m(v) = u} e^{g_m(v)}\, K h(v)\\
& =    & e^{ C' \, \Lipf(g)}\,C'K \Lip(g) \, d(u,u')\, (L_g^m h)(u)\\
& =    &  e^{ C' \, \Lipf(g)}\,C'K \Lip(g) \, d(u,u')\, \lambda^m h(u)\\
& \leq & \lambda^m\, C'K^2 e^{ C' \, \Lipf(g)}\, \Lip(g) \, d(u,u') .
\end{eqnarray*}
Thus, for every integer $m$ the function $\frac{1}{\lambda^m} L_g^m 1 \in \clip(\hU)$ and 
$\frac{1}{\lambda^m} \Lip(L_g^m 1) \leq C'K^2 e^{ C' \, \Lipf(g)}\, \Lip(g)$. As mentioned above this proves that
the eigenfunction $h \in \clip(\hU)$. 

Using this with $g = \ft$ proves that $h_{at} \in \clip(\hU)$ for all $|a| \leq a_0$ and $t \geq 1/a_0$.
However the above estimate for $\Lip(h_{at})$ would be of the form $\leq C \, e^{C \, t} \, t$ for some constant $C > 0$, 
which is not good enough.

We will now show that, taking $a_0 > 0$ sufficiently small, we have $\Lip(h_{at}) \leq C t$ for some constant
$C > 0$ independent of $a$ and $t$.

Using (3.2) and choosing $a_0 > 0$ sufficiently small, we have $\lambda_{at} \gamma > \hgamma$
for all $|a| \leq a_0$ and $t > 1/a_0$. Fix an integer $m_0 \geq 1$ so large that
$\frac{C_0^2}{c_0 \hgamma^{m}} < \frac{1}{2}$ for $m \geq m_0$. There exists a constant $d_0 > 0$ depending on
$m_0$ such that for any $u,u'$ belonging to the same $U_i$ but not to the same cylinder of length $m_0$ we have
$d(u,u') \geq d_0$. For such $u,u'$ we have 
$$\frac{|h_{at}(u) - h_{at}(u')|}{d(u,u')} \leq \frac{2\|h_{at}\|_0}{d_0} \leq \frac{2C_0}{d_0} .$$
So, to estimate $\Lip(h_{at})$ it is enough to consider pairs $u,u'$ that belong to the same cylinder of length $m_0$.

Fix for a moment $a,t$ with $|a| \leq a_0$ and $t \geq 1/a_0$. Set 
$$L = \sup_{u\neq u'}\frac{|h_{at}(u) - h_{at}(u')|}{d(u,u')} ,$$ 
where the supremum is taken over all pairs $u \neq u'$ that belong to the same cylinder of length $m_0$. 
If $L < \Lip(h_{at})$, then the above implies
$$\Lip(h_{at}) \leq \frac{2C_0}{d_0} \leq \frac{2C_0}{d_0}\, t .$$

Assume that $L = \Lip(h_{at})$. 
Then there exist $u, u'$ belonging to the same cylinder of length $m_0$ such that
\be
\frac{3L}{4} < \frac{|h_{at}(u) - h_{at}(u')|}{d(u,u')} .
\ee

Fix such a pair $u,u'$. Let  $m \geq m_0$ be an integer. For any $v\in \hU$ with
$\sigma^m(v) = u$, denote by $v' = v'(v)$ the unique $v'\in \hU$ in the cylinder of length $m$ containing $v$ such
that $\sigma^m(v') = u'$. By (1.1),
$$d (\sigma^j(v),\sigma^j(v')) \leq \frac{1}{c_0\, \gamma^{m-j}}\, d (u,u') \quad, \quad j = 0,1, \ldots, m-1 $$
so
$$|\ft^m(v) - \ft^m(v')| \leq \sum_{j=0}^{m-1} |\ft(\sigma^j(v)) - \ft(\sigma^j(v'))|
\leq \Con \Lip(\ft) \, d(u,u') \leq \Con \, t\, d(u,u') .$$ 
At the same time, by property (i), $\|\ft\|_0\leq T''$ for some constant $T'' > 0$, so
$$|\ft^m(v) - \ft^m(v'(v))| \leq 2m \|\ft\|_0 \leq 2m T'' .$$
Similarly, 
$$|(P+a) \tau^m(v) - (P+a)\tau^m(v')| \leq \Con \,  d(u,u') \leq T'' ,$$ 
assuming $T'' > 0$ is chosen sufficiently large. Thus,
\begin{eqnarray*}
&      & \left|e^{(\ft - (P+a)\tau)^m(v') - (\ft - (P+a)\tau)^m(v)} -1\right|\\
& \leq & e^{3m T''}\, \left|(\ft - (P+a)\tau)^m(v) - (\ft - (P+a)\tau)^m(v')\right|
\leq e^{3m T''}\, \Con \, t\, d(u,u') .
\end{eqnarray*}

Using  $L^m_{\ft - (P+a)\tau}h_{at} = \lambda_{at}^m h_{at}$, we obtain
\begin{eqnarray*}
&        &\lambda_{at}^m\, |h_{at} (u) - h_{at}(u')| 
 =  \left| \sum_{\sigma^m v = u} e^{(\ft - (P+a)\tau)^m(v)}\, h_{at}(v) -  
\sum_{\sigma^m v = u} e^{(\ft - (P+a)\tau)^m(v'(v))}\, h_{at}(v') \right| \\
& \leq & \sum_{\sigma^m v = u} e^{(\ft - (P+a)\tau)^m(v)}\, |h_{at}(v) -  h_{at}(v')| 
       + \|h_{at}\|_0 \sum_{\sigma^m v = u} \left|e^{(\ft - (P+a)\tau)^m(v)} - e^{(\ft - (P+a)\tau)^m(v')}\right| \\
& \leq & \frac{\Lip(h_{at})\, d(u,u')}{c_0\gamma^m}\sum_{\sigma^m v = u} e^{(\ft - (P+a)\tau)^m(v)}\\
&      & + C_0\, \sum_{\sigma^m v = u} e^{(\ft - (P+a)\tau)^m(v)} 
\left|1 - e^{(\ft - (P+a)\tau)^m(v') - (\ft - (P+a)\tau)^m(v)}\right| \\
& \leq & \frac{L \, d(u,u')}{c_0\gamma^m}\sum_{\sigma^m v = u}  e^{(\ft - (P+a)\tau)^m(v)}
       +  C_0 e^{3m T''}\, \Con \, t\, d(u,u') \sum_{\sigma^m v = u} e^{(\ft - (P+a)\tau)^m(v)} \\
& \leq & \left( \frac{L}{c_0\gamma^m} + C_0 e^{3m T''}\, \Con \, t  \right)\, d(u,u')\, 
\sum_{\sigma^m v = u}  e^{(\ft - (P+a)\tau)^m(v)}\, C_0 h_{at}(v)\\
& =    &  \left( \frac{L}{c_0\gamma^m} + C_0 e^{3m T''}\, \Con \, t  \right)\, d(u,u')\, C_0 \lambda_{at}^m h_{at}(u)
\leq  \left( \frac{L}{c_0\gamma^m} + C_0 e^{3m T''}\, \Con \, t  \right)\, d(u,u')\, C^2_0 \lambda_{at}^m .
\end{eqnarray*}
This, (8.1) and the choice of $m_0$ imply
$$\frac{3L}{4} < \frac{L C_0^2}{c_0\gamma^m} + C^3_0 e^{3m T''}\, \Con \, t 
\leq \frac{L}{2} + C^3_0 e^{3m T''}\, \Con \, t .$$
This is true for all $m \geq m_0$. In particular for $m = m_0$ we get
$$\frac{L}{4} <  C^3_0 e^{3m_0 T''}\, \Con \, t ,$$
and so $\Lip(h_{at}) = L \leq \Con \; t$.
\endofproof

\bs

\noindent
{\it Proof of Lemma} 5. 
(a) Let  $u, u' \in \hU_i $ for some $i = 1, \ldots,k$ and let $m \geq 1$ be an integer. For any $v\in \hU$ with
$\sigma^m(v) = u$, denote by $v' = v'(v)$ the unique $v'\in \hU$ in the cylinder of length $m$ containing $v$ such
that $\sigma^m(v') = u'$. 
Then
\begin{eqnarray}
|\fat^m(v) - \fat^m(v')|  \leq  \sum_{j=0}^{m-1} |\fat(\sigma^j(v)) - \fat (\sigma^j(v'))| 
 \leq  \frac{ T t}{c_0\, (\gamma-1)}\, d (u,u') \leq C_1\,t\,  D (u,u') 
\end{eqnarray}
for some constant $C_1 > 0$. Similarly, 
\be
|\gt^m(v) - \gt^m(v')| \leq C_1 \, t\, D(u,u') .
\ee

Also notice that 
if $D(u,u') = \diam(\cc')$ for some cylinder $\cc' = C[i_{m+1}, \ldots, i_p]$,then $v,v'(v) \in \cc'' = C[i_0,i_1, \ldots,i_p]$ 
for some cylinder $\cc''$ with $\sigma^m(\cc'') = \cc'$, so 
$$D(v,v') \leq \diam(\cc'')  \leq \frac{1}{c_0\, \gamma^m} \, \diam(\cc')  = \frac{D(u,u')}{c_0\, \gamma^m} .$$

Using the above, $\diam (U_i) \leq 1$, the definition of $\matc$, we get
\begin{eqnarray*}
&        &\frac{|(\matc^m H)(u) - (\matc^m H)(u')|}{\matc^m H (u')} 
 =      \frac{\di \left| \sum_{\sigma^m v = u} e^{\fat^m(v) + c \gt^m(v)}\, H(v) -  
\sum_{\sigma^m v = u} e^{\fat^m(v') + c \gt^m(v')}\, H(v') \right|}{\matc^m H (u')} \\
& \leq & \frac{\di \left| \sum_{\sigma^m v = u} e^{\fat^m(v)+ c \gt^m(v)}\, (H(v) -  H(v'))\right|}{\matc^m H (u')}   +
\frac{\di \sum_{\sigma^m v = u}  \left|e^{\fat^m(v) + c \gt^m(v)}-  e^{\fat^m(v') + c \gt^m(v')}\right| \, H(v')}{\matc^m H (u')} \\
& \leq & \frac{\di  \sum_{\sigma^m v = u} e^{\fat^m (v) + c \gt^m(v)}\, E\, H(v')\, D (v,v')}{\matc^m H (u')} \\
&      &  + \frac{\di \sum_{\sigma^m v = u} \left| e^{[\fat^m(v) + c \gt^m(v)] - [\fat^m(v')+ c \gt^m(v')]} -1 \right|  \,
e^{\fat^m(v') + c \gt^m(v')}\, H(v')}{\matc^m H (u')} .
\end{eqnarray*}
Using (8.2) and (8.3) and assuming $\eta_0 \leq 1$, one obtains
\be
|\fat^m(v) + c \gt^m(v)] - [\fat^m(v')+ c \gt^m(v')| \leq 2C_1 t\, D(u,u') \leq 2C_1 t ,
\ee
and therefore
$$\left| e^{[\fat^m(v) + c \gt^m(v)] - [\fat^m(v')+ c \gt^m(v')]} -1 \right| \leq e^{2C_1t} 2C_1 t\, D(u,u') .$$
However (8.4) is not good enough to estimate the first term in the right-hand-side above. Instead we use (3.3) and (3.4) to get
\begin{eqnarray}
&      & |\fat^m(v) + c \gt^m(v)] - [\fat^m(v')+ c \gt^m(v')|\nonumber\\
& \leq & |\ft^m(v) - \ft^m(v)| + |P-a| \, |\tau^m(v) - \tau^m(v')| 
+ |(h_{at}(v) - h_{at}(u)) - (h_{at}(v') - h_{at}(u')|\nonumber\\
&      & + a_0 |\gt^m(v) - \gt^m(v')|\nonumber\\
& \leq & 2m \|\ft - \f0\|_0 + |\f0^m(v) - \f0^m(v')| + \Con \, D(u,u') + 4C_0 + 2m a_0 \|\gt- g\|_0\nonumber\\
& \leq & \Con \, D(u,u') + C_2 m a_0 \leq C_2 + C_2 m\, a_0 
\end{eqnarray}
for some constant $C_2 > 0$. We will now assume that $a_0 > 0$ is chosen so small that
\be
e^{C_2a_0} < \gamma/\hgamma .
\ee
Hence
\begin{eqnarray*}
&      & \frac{|(\matc^m H)(u) - (\matc^m H)(u')|}{\matc^m H (u')} \\
& \leq & \frac{E\, D (u,u')}{c_0\gamma^m}\, \frac{\di  \sum_{\sigma^m v = u}  
e^{[\fat^m(v)+ c \gt^m(v)] - [\fat^m(v')+ c \gt^m(v')]} e^{\fat^m(v')+ c \gt^m(v')} \,  H(v')}{\matc^m H (u')}\\
&        & \:\:\:  + e^{2C_1 t}\, \frac{\di \sum_{\sigma^m v = u} 2C_1 t\,e^{\fat^m(v'(v))}\, H(v'(v))}{\matc^m H (u')} \\
& \leq &  e^{C_2}\,e^{C_2m a_0 } \frac{E\, D (u,u')}{c_0\gamma^m} + 2C_1 t e^{2C_1 t}\, D (u, u')
 \leq  A_0 \, \left[ \frac{E}{\hgamma^m} + e^{A_0 t}\, t \right]\, D (u,u') ,
\end{eqnarray*}
for some constant $A_0 > 0$ independent of $a$, $c$, $t$, $m$ and $E$.

\ms

(b) Let $m \geq 1$ be an integer and $u,u'\in \hU_i$ for some $i = 1, \ldots,k$. 
Using the notation $v' = v'(v)$ and the constant $C_2 > 0$  from part (a) above, where $\sigma^m v = u$ and
$\sigma^mv' = u'$, and some of the estimates from the proof of part (a), we get
\begin{eqnarray*}
&      & |\labtz^m h(u) - \labtz^m h(u')|\\
&  =   &    \di \left| \sum_{\sigma^m v = u} \left( e^{\fat^m(v) + c \gt^m(v) - \i b \tau^m(v) + \i w \gt^m(v)}\, h(v) -
  e^{\fat^m(v') + c \gt^m(v') - \i b \tau^m(v') + \i w \gt^m(v')}\, h(v')\right)\right|\\
& \leq & \di \left| \sum_{\sigma^m v = u} e^{\fat^m(v) + c \gt^m(v) - \i b \tau^m(v) + \i w \gt^m(v)}\, [h(v) - h(v')]\right|\\
&      &   \di  +  \sum_{\sigma^m v = u}  \left| e^{\fat^m(v) + c \gt^m(v)} -
  e^{\fat^m(v') + c \gt^m(v')}\right|\, |h(v')|\\
&      &   \di  +   \sum_{\sigma^m v = u} \left| e^{ - \i b \tau^m(v) + \i w \gt^m(v)} -
  e^{- \i b \tau^m(v') - \i w \gt^m(v')}\right|\, e^{\fat^m(v') + c \gt^m(v')} |h(v')|\\
& \leq & \di \sum_{\sigma^m v = u} e^{\fat^m(v) + c \gt^m(v)}\, | h(v) - h(v')|\\
&  +   & \sum_{\sigma^m v = u} \left| e^{[\fat^m(v) + c \gt^m(v)] 
- [\fat^m(v') + c \gt^m(v')]} - 1\right|\, e^{\fat^m(v') + c \gt^m(v') }\, |h(v')| \\
&      &  \sum_{\sigma^m v = u} \left(|b|\, |\tau^m(v) - \tau^m(v')| + |w|\,  |\gt^m(v) -\gt^m(v')|\right)\, 
e^{\fat^m(v') + c \gt^m(v')} |h(v')|\\
\end{eqnarray*}

Using  the constants $C_1, C_2  > 0$ from the proof of part (a), (8.5) and (8.6) we get
\begin{eqnarray*}
        \sum_{\sigma^m v = u} e^{\fat^m(v) + c \gt^m(v)}\, | h(v) - h(v')|
& \leq & e^{C_2}\,e^{C_2m a_0}\frac{E\, D (u,u')}{c_0\gamma^m}\sum_{\sigma^m v = u} e^{\fat^m(v') + c \gt^m(v')}\, H(v')\\
& \leq & \frac{e^{C_2} E}{c_0 \hgamma^m} (\matc^m H)(u') \, D(u,u') .
\end{eqnarray*}
This, (8.3) and (8.5) imply
\begin{eqnarray*}
&      & |\labtz^m h(u) - \labtz^m h(u')|\\
& \leq & \frac{e^{C_2} E}{c_0 \hgamma^m} (\matc^m H)(u') \, D(u,u') + e^{2C_1t} 2C_1 t\, D(u,u')\, (\matc^m |h|)(u')
+ \left( \Con \, |b| + |w|C_1 \, t \right) D(u,u')
\end{eqnarray*}
Thus, taking the constant $A_0 > 0$ sufficiently large we get
$$|(\labtz^N h)(u) - (\labtz^N  h)(u')|\\
 \leq  A_0  \left(\frac{E}{\hgamma^m} (\matc^m H)(u') + ( |b| + e^{A_0t}t + t|w|) (\matc^m |h|)(u')\right)\, D(u,u') ,$$
which proves the assertion.
\endofproof

\bs

As in \cite{D} and \cite{St2} we need the following lemma whose proof is omitted here, since it is very similar to the 
proof of Lemma 5 given above.

\ms

\begin{lem}
Let $\beta \in (0,\alpha)$. There exists a constants $A'_0 > 0$ such that for all $a,b,c, t,w\in \R$ with 
$|a|, |c|,  1/|b|, 1/t\leq a_0$ such that {\rm (4.1)} hold, and all positive 
integers $m$  and all  $h \in C^\beta (U)$ we have
$$|\labtz^m h(u) - \labtz^m h(u')| \leq A'_0\left[ \frac{|h|_\beta}{\hgamma^{m\beta}} 
+ |b| \, (\matc^m |h| )(u')\right]\, (d (u,u'))^\beta$$
for all $u,u'\in U_i$. 
\end{lem}

\ms

We will derive Theorem 5(b) from Theorem 5(a), proved in Sect. 4,  and Lemma 17 above.

\bs

\noindent
{\it Proof of Theorem} 5(b). We essentially repeat the proofs
of Corollaries 2 and 3 in \cite{D} (cf. also Sect. 3 in \cite{St1}).

Let $\ep > 0$, $B > 0$ and $\beta \in (0,\alpha)$. Take $\hrho\in (0,1)$, $a_0 > 0$, $b_0 > 0$, $A_0 > 0$ and $N$  
as in Theorem 2(a). We will assume that $\hrho \geq \frac{1}{\gamma_0}$. Let $a, b,c,w\in \R$  be such that $|a|, |c| \leq a_0$ 
and $|b| \geq b_0$. Let $t > 0$ be such that $1/t^{\alpha - \beta} \leq a_0$. Assume that (4.1) hold and set $z = c + \i w$.

First, as in \cite{D} (see also Sect. 3 in \cite{St1}) one derives from Theorem 5(a) and Lemma 17 (approximating functions
$h \in C^\beta(\hU)$ by Lipschitz functions as in Sect. 3)
that there exist constants $C_3 > 0$ and $\rho_1 \in (0,1)$ such that
\be
\|\labtz^{n} h\|_{\beta,b} \leq C_{3} |b|^\ep \rho_{1}^n \quad , \quad n \geq 0 ,
\ee
for all $h \in C^\beta(\hU)$.

Next, given $h \in C^\beta(\hU)$, we have
$$\labtz^n(h/h_{at}) = \frac{1}{\lambda^n_{at} \, h_{at} } \, L_{\ft - (P + a+\i b)\tau + z\gt} h ,$$
so by (8.7) we get
\begin{eqnarray*}
\|L^n_{\ft - (P + a+\i b)\tau + z \gt}h\|_{\beta,b} 
& \leq & \lambda_{at}^n \|h_{at}\, \labtz^n(h/h_{at})\|_{\beta,b}\\
& \leq & \Con (\lambda_{at} \rho_1)^n |b|^\ep\, \|h/h_{at}\|_{\beta,b}
\leq  \Con \, \rho_2^n \, |b|^\ep\, \|h\|_{\beta,b}\;,
\end{eqnarray*}
where $\lambda_{at}\rho_1 \leq e^{2C_0a_0}\rho_2 = \rho_2 < 1$, provided $a_0 > 0$ is small enough.

We will now approximate $L_{f- (P+a+\i b)\tau + z g}$ by $L_{\ft - (P + a+\i b)\tau +\gt}$
in two steps. First, using the above it follows that
\begin{eqnarray*}
\|L^n_{f - (P + a + \i b) \tau + c g + \i w \gt} h\|_{\beta,b} 
& =    & \left\|L^n_{\ft - (P + a + \i b)\tau + z\gt} \left(e^{(f^n-\ft^n) + c (g^n - \gt^n)} h\right)\right\|_{\beta,b}\\
& \leq & \Con \, \rho_2^n \, |b|^\ep\,  \left\|e^{(f^n-\ft^n) + c (g^n -\gt^n)} h \right\|_{\beta,b} .
\end{eqnarray*}
Choosing the constant $C_4 > 0$ appropriately, $\|f-\ft\|_0 \leq C_4 \, a_0$ and 
$|f - \ft|_\beta \leq C_4 /t^{\alpha-\beta} \leq C_4 a_0$, so 
$\|f^n - \ft^n\|_0 \leq n\, \|f-\ft|_0 \leq C_4 n a_0$, and similarly
$|f^n - \ft^n|_\beta \leq C_4 n a_0$. Similar estimates hold for $g^n - \gt^n$.  Thus,
$$\|e^{(f^n-\ft^n) + c (g^n - \gt^n)} h \|_0 \leq e^{C_4 n a_0}\|h\|_0$$
and
\begin{eqnarray*}
|e^{(f^n-\ft^n)  + c (g^n -\gt^n)} h |_\beta 
& \leq &  \|e^{(f^n-\ft^n) + c (g^n -\gt^n)}\|_0\, |h|_\beta + |e^{(f^n-\ft^n)  + c(g^n -\gt^n)} |_\beta \, \|h\|_\infty\\
& \leq & e^{C_4 n a_0} |h|_\beta  + e^{C_4 n a_0}\, |(f^n-\ft^n) + c (g^n -\gt^n)|_\beta\, \|h\|_\infty\\
& \leq & C'_5 \, n \,e^{C_4 n a_0}\,\|h\|_\beta .
\end{eqnarray*}
Combining this with the previous estimate gives
$$\| e^{(f^n-\ft^n) + c(g^n -\gt^n)} h \|_{\beta,b} \leq C''_5 \, n \,e^{C_4 n a_0}\,\|h\|_\beta ,$$
so 
$$\|L^n_{f - (P + a + \i b) \tau + cg + \i w \gt} h\|_{\beta,b} \leq
C_5 \, \rho_2^n \, |b|^\ep\, n \,e^{C_4 n a_0}\,\|h\|_{\beta,b} .$$
Taking $a_0 > 0$ sufficiently small, we may assume that $\rho_2\, e^{C_4 a_0} < 1$. Now take an arbitrary $\rho_3$ with
$\rho_2\, e^{C_4a_0} < \rho_3 < 1$. Then we can take the constant $C_6 > 0$
so large that $n \,\rho^n_2\, e^{C_4 n a_0} \leq C_6 \rho_3^n$ for all  integers $n \geq 1$. This gives
\be
\|L^n_{f - (P + a + \i b) \tau + cg + \i w \gt} h\|_{\beta,b} \leq
C_6 \, \rho_3^n \, |b|^\ep\, \|h\|_{\beta,b} \quad , \quad n \geq 0 .
\ee

Using the latter we can write
\begin{eqnarray*}
\|L^n_{f - (P + a + \i b )\tau + z g} h\|_{\beta,b} 
& =    & \left\|L^n_{f - (P + a + \i b) \tau + c g + \i w \gt} \left(e^{\i w (g^n -\gt^n)} h\right)\right\|_{\beta,b}\\
& \leq & C_6 \, \rho_3^n \, |b|^\ep\,  \left\|e^{\i w (g^n -\gt^n)} h \right\|_{\beta,b} .
\end{eqnarray*}
However, $\|e^{\i w (g^n -\gt^n)} h \|_0 = \|h\|_0$, $|g-\gt|_\beta \leq C_4/ t^{\alpha - \beta} \leq C_4 a_0 \leq 1$
(assuming $a_0 > 0$ is sufficiently small),  and by (4.1), $|w| \leq B |b|^\mu \leq B |b|$, so
\begin{eqnarray*}
|e^{\i w (g^n -\gt^n)} h |_\beta 
& \leq &  \|e^{\i w (g^n -\gt^n)}\|_0\, |h|_\beta + |e^{\i w (g^n -\gt^n)} |_\beta \, \|h\|_\infty\\
& \leq &  |h|_\beta + |w|\, |g^n -\gt^n|_\beta\, \|h\|_\infty\\
& \leq & \|h\|_\beta + B n |b| \|h\|_\infty .
\end{eqnarray*}
Thus,
$$\|e^{\i w(g^n -\gt^n)} h \|_{\beta,b} 
= \|e^{\i w(g^n -\gt^n)} h \|_0 + \frac{1}{|b|} |e^{\i w(g^n -\gt^n)} h |_\beta  \leq 2 B n \|h\|_{\beta,b} ,$$
and therefore
$$\|L^n_{f - (P + a + \i b )\tau + z g}h\|_{\beta,b} \leq C_7\, \rho_3 ^n \, |b|^\ep\, n\, \|h\|_{\beta,b} .$$
Now taking an arbitrary $\rho$ with $\rho_3 < \rho < 1$ and taking the constant $C_8 > C_7$
sufficiently large, we get
$$\|L^n_{f - (P + a + \i b )\tau + z g}h\|_{\beta,b}  \leq C_8 \, \rho^n \, |b|^\ep\, \|h\|_{\beta,b}$$
for all integers $n \geq 0$.
\endofproof

\bs

{\footnotesize


\begin{thebibliography}{99}

\bibitem{B1} R. Bowen, Equilibrium states and the ergodic theory of Anosov diffeomorphisms, Lect. Notes in Maths.
{\bf 470}, Springer-Verlag, Berlin, 1975.

\bibitem{B2} R. Bowen, {\em Symbolic dynamics for hyperbolic flows}, Amer. J. Math. {\bf 95} (1973), 429-460.

\bibitem{BR} R. Bowen and D. Ruelle, {\em  The ergodic theory of Axiom A flows}, Invent. Math. {\bf 29} (1975), 181-202.

\bibitem{D} D. Dolgopyat, {\em Decay of correlations in Anosov flows}, Ann. Math. {\bf 147} (1998), 357-390.

\bibitem{HO} J. M. Hannay and A. M. Ozorio de Almeida, {\em Periodic orbits and a correlation function for the 
semiclassical density of states}, J. Phys. A {\bf 17} (1984), 3429-3440.

\bibitem{KH} A. Katok and  B. Hasselblatt, {\em Introduction to the Modern Theory of
Dynamical Systems}, Cambridge Univ. Press, Cambridge 1995.

\bibitem{KS} A. Katsuda and T. Sunada, {\em Closed orbits in homology class}, Publ. math\'ematiques d'IHES, {\bf 71} (1990), 5-32.

\bibitem{La} S. Lalley, {\em Distribution of period orbits of symbolic and Axiom A flows}, Adv. Appl. Math. {\bf 8} (1987), 154-193.

\bibitem{N} F. Naud, {\em Expanding maps on Cantor sets and analytic continuation of zeta function},  
Ann. Sci. Ec. Norm. Sup. {\bf 38} (2005), 116-153.

\bibitem{Pa} W. Parry, {\em Synchronization of canonical measures for hyperbolic attractors}, Comm. Math. Phys. 
{\bf 106} (1986), 267-275.

\bibitem{PP} W. Parry and M. Pollicott,  {\em Zeta functions and  the periodic orbit
structure of hyperbolic dynamics}, Ast\'erisque {\bf 187-188}, (1990).

\bibitem{PeS} V. Petkov and L. Stoyanov, {\em Sharp large deviations for some hyperbolic systems},
Ergod. Th. \& Dyn. Sys., doi: 10.1017/etds.2013.48.

\bibitem{Po} M. Pollicott, {\em On the rate of mixing of Axiom A flows}, Invent. Math. {\bf 81} (1985), 413-426.

\bibitem{Po1} M. Pollicott, {\em A note on exponential mixing for Gibbs measures and counting weighted periodic orbits for 
geodesic flows}, Preprint 2014.

\bibitem{PoS1} M. Pollicott and R. Sharp, {\em Exponential error terms for growth functions on negatively curved surfaces}, 
Amer. J. Math. {\bf 120} (1998), 1019-1042.

\bibitem{PoS2} M. Pollicott and R. Sharp, {\em Large deviations, fluctuations and 
shrinking intervals},  Comm. Math. Phys. {\bf 290} (2009), 321-324.

\bibitem{PoS3} M. Pollicott and R. Sharp, {\em On the Hannay-Ozorio de Almeida sum formula}, Dynamics, games and science. II, 575-590, Springer Proc. Math., 2, Springer, Heidelberg, 2011.

\bibitem{R} D. Ruelle, {\em An extension of the theory of Fredholm determinants}, Publ. Math. IHES, {\bf 72} (1990), 175-193.

\bibitem{St1} L. Stoyanov, {\em Spectrum of the Ruelle operator and exponential decay of 
correlations for open billiard flows}, Amer. J.  Math. {\bf 123}, (2001), 715-759.

\bibitem{St2} L. Stoyanov, {\em Spectra of Ruelle transfer operators for Axiom A flows}, 
Nonlinearity, {\bf 24} (2011), 1089-1120.

\bibitem{St3} {L. Stoyanov}, {\em Pinching conditions, linearization and regularity of  
Axiom A flows}, {\em Discr.  Cont. Dyn. Sys. A}, {\bf 33} (2013), 391-412.



\bibitem{Wa} S. Waddington, {\em Large deviations for Anosov flows}, Ann. Inst. H. Poincar\'e, Analyse non-lin\'eaire, 
{\bf 13}, (1996), 445-484.

\bibitem{W} P. Wright, {\em Ruelle's lemma and Ruelle zeta functions}, Asymptotic Analysis, {\bf 80} (2012), 223-236.



\end{thebibliography}
 \end{document}